\documentclass[12pt]{amsart}
\usepackage{amsmath,amssymb,amsfonts}
\usepackage{amsthm}
\usepackage{array}
\usepackage{resizegather}
\usepackage[shortlabels]{enumitem}
\usepackage{color}
\usepackage{graphicx}

\theoremstyle{plain}

\theoremstyle{definition}

\numberwithin{equation}{section}

\newcommand{\scr}[1]{\mbox{\scriptsize $#1$}}

\newcommand{\La}[1]{\mbox{\Large $#1$}}
\newcommand{\LA}[1]{\mbox{\LARGE $#1$}}

\begin{document}

\title[Grassman manifolds as subsets of Euclidean spaces]
{Grassman manifolds as subsets of Euclidean spaces}

\author[A. Machado, I. Salavessa]{
Armando Machado$^{\dag}$, Isabel Salavessa$^{\ddag}$\\
 \tiny{ \em Faculty of Sciences, University of Lisboa \&
Centro de Matem\'{a}tica e Aplica\c{c}\~{o}es Fundamentais, Av.\  Prof.\ Gama Pinto 2, 1049-003 Lisboa, Portugal \em }
}
\address{}
\email{}
\address{}
\email{}
\footnote{\noindent $\dag$ email: ahmachado@ciencias.ulisboa.pt\\
$\ddag$ current address: Centro de F\'{\i}sica e Engenharia de Materiais 
Avan\c{c}ados, Edif\'{\i}cio Ci\^{e}ncia, Piso 3,
Instituto Superior T\'{e}cnico, Universidade de Lisboa, Av. Rovisco Pais, 1049-001 Lisboa, Portugal;\\ email: isabel.salavessa@tecnico.ulisboa.pt }

\date{1984}

\maketitle

\section{Introduction}
Let $E$ be a Euclidean space. Following Palais, we identify each vector subspace F of E with
 the orthogonal projection $\pi_F:E\to F$. In this way, the Grassman manifold $G(E)$ of all
 vector subspaces of E appears as a submanifold of the Euclidean space $L(E;E)$ of all
 linear maps from $E$ into $E$ (with the Hilbert-Schmidt inner product). The aim of this
 paper is to present some explicit formulas concerning the differential geometry of $G(E)$
 as a submanifold of $L(E;E)$. Most of these formulas extend naturally to the case where $E$
 is an infinite dimensional Hilbert space, although in this case there is no natural inner
 product in $L(E;E)$.

\section{Notation and Preliminaries }

Let $E$ and $F$ be finite or infinite dimensional Hilbert spaces. We will denote by $L(E;F)$
 the vector space of all continuous linear maps from $E$ into $F$. If $\xi\in L(E;F)$, we
will denote by $\xi^*\in L(F;E)$ its \em adjoint \em  linear map, the one defined by the
identity 
\[\langle\xi(x),y\rangle \,=\, \langle x,\xi^*(y)\rangle.\] 
The following identities will be used quite often:
\begin{equation}
\xi^{**}\,=\,\xi;\quad (\eta\circ\xi)^*\,=\, \xi^*\circ \eta^*; \quad id_E^*\,=\, id_E.
\end{equation}
A linear map $\xi\in L(E;E)$ is \em self-adjoint \em  if $\xi^*=\xi$. The map $L(E;F)\to
 L(F;E)$, $\xi\to \xi^*$, is a \em real \em linear map (even if $E$ and $F$ are complex
 spaces, it is not a complex linear map) and the set $L_{sa}(E;E)$ of self-adjoint linear
 maps is a real vector subspace of $L(E;F)$. 

In case $E$ or $F$ is infinite dimensional, we will look on  $L(E;F)$ merely as a 
Banach space (with the sup norm). In case $E$ and $F$ are finite dimensional, we take 
in the finite dimensional vector space $L(E;F)$ the \em
 Hilbert-Schmidt inner product, \em defined by 
\begin{equation}
\langle \xi,\eta\rangle \, =\, \sum_{1\leq k\leq n}\langle \xi(x_k),\eta(x_k)\rangle,
\end{equation}
where $x_1, \ldots, x_k$ is an arbitrary orthonormal basis of $E$. We will use the
  following identities concerning these inner products,
\begin{equation}
\langle \xi, \eta\rangle \, =\, \langle \eta^*, \xi^*\rangle; \quad
\langle \lambda, \mu^*\circ \eta\rangle \, =\, \langle \mu\circ \lambda, \eta\rangle \,
=\, \langle \mu, \eta \circ \lambda^*\rangle.
\end{equation}

The word ``manifold'' will always mean an embedded submanifold of some finite dimensional or
 Banach vector space $B$ and the tangent vector spaces will be considered as vector
 subspaces of the ambient vector space $B$. In fact, one can even define, for each point $a$
 of an arbitrary subset $M$ of $B$, a notion of tangent vector subspace $T_a(M)$, which
 behaves well with respect to differentiability (see, for example, [2]). In the same spirit,
 by vector bundle we will mean a vector sub-bundle of a constant one. A vector bundle 
       $\underline{E}$
 with basis $M$ will be a family $(E_x)_{x\in M}$, where each $E_x$ is a vector subspace of a
 fixed finite dimensional or Banach vector space $E$, verifying the usual properties, and we
 will use the same symbol $\underline{E}$ to denote the corresponding subset of $M \times E$.
 It will be useful to allow a vector bundle to have as basis an arbitrary subset M of a
 finite dimensional or Banach vector space B. 

 If  $\underline{E} = (E_x)_{x\in M}$ is a
 vector bundle with $E_x\subset E$, we identify a \em connection \em in $\underline{E}$ by
 its second fundamental form at each point $x\in M$, which is a bilinear map 
$\theta_x:E_x\times T_x(M) \to E$  such that
\begin{equation}
(u,\theta_x(w,u)) \in T_{(x,w)}(\underline{E}).
\end{equation}
For each smooth section $W =(W_x)_{x\in M}$ of $\underline{E}$, the covariant derivative
$\nabla W_x(u)$ is  given by the formula
\begin{equation}
 \nabla W_x(u)= DW_x(u) - \theta_x(W_x,u). 
\end{equation}
If $E$ is a Hilbert space, the \em metric connection \em of $\underline{E}$ is the one
 defined by the condition that $\theta_x(w,u)$ is orthogonal to the fibre $E_x$; if 
$\pi_x:E\to E_x$ is the orthogonal projection, then $x\to \pi_x$ is a smooth map from $M$
 into $L(E;E)$ and we have the following formula for this connection,
\begin{equation}
\theta_x(w,u)=D\pi_x(u)(w).
\end{equation}
 
We will use also the following characterization of the curvature tensor of a connection 
$\theta$  in the vector bundle $\underline{E}= (E_x)_{x\in M}$, where $M \subset B$ is a
 manifold and $E_x \subset E$: assuming that $x\to \hat{\theta}_x$ is a smooth map from $M$
 into the space $L(E,B;E)$ of bilinear maps, such that each $\theta_x$ is a restriction of 
$\hat{\theta}_x$, the curvature tensor is the trilinear map 
\[  R_x: T_x(M) \times T_x(M) \times E_x \to  E_x\]
defined by
\begin{equation}
         R_x(u,v,w) = D\hat{\theta}_x(u)(w,v)-D\hat{\theta}_x(v)(w,u) +
 \hat{\theta}_x(\theta_x(w,u),v)- \hat{\theta}_x(\theta_x(w,v),u).
\end{equation}

\section{The Grassman Manifolds}

Let $E$ be a finite or infinite dimensional real Hilbert space. For each closed 
vector subspace $F \subset E$, we will denote by $\pi_F$ the orthogonal projection from 
$E$ onto $F$. We have hence a natural bijective map between the set of closed 
vector subspaces of $E$ and the set of orthogonal projections. We will denote 
by $G(E)$ the subset of $L(E;E)$ whose elements are the orthogonal projections 
onto closed subspaces, and we will call $G(E)$ the \em Grassman manifold \em  of $E$.
The fact that $G(E)$ is indeed a manifold is proved in Akin [1], who attributes 
this result to Palais (unpublished preprint), but we will sketch here an 
independent proof.

    The following characterization of the elements of $G(E)$ is well known: 

\subsection{} A linear map $\xi \in L(E;E)$ belongs to $G(E)$ if and only if it is 
self-adjoint and verifies $\xi\circ \xi = \xi$.

    We can consider a morphism from the constant vector bundle $E_{G(E)}$, with
basis $G(E)$ and fibre $E$, into itself, associating to each $\xi \in G(E)$ the linear 
map $\xi:E \to E$. The fact that the image of an idempotent morphism is a vector 
bundle allows us to state:

\subsection{} There exists a \em tautological vector bundle \em with basis $G(E)$, whose
 fibre in each $\pi_F$ is $F$.
 
Using formula (2.6) for the metric connection, we deduce:

\subsection{} The metric connection of the tautological vector bundle is defined by
\[ \theta_{\xi}(w,n) \,=\, \eta(w),\]
for each $\xi\in G(E)$, $w\in \xi(E)$  and $\eta\in T_{\xi}(G(E))$.\\

        As a corollary of the local constancy of the dimension of the fibres of a 
vector bundle, we see that, for each $n$, the subset $G_n(E)$ of $G(E)$, whose 
elements are the $\pi_F$ such that $F$ is $n$-dimensional, is open in $G(E)$.

        Let $F\subset E$ be a fixed closed vector subspace. It is a well known simple 
linear algebra result that, for each closed vector subspace $G\subset E$, the 
following two properties are equivalent:
\[\begin{array}{rll}
        (a) &E=F^{\bot}\oplus G  \quad\mbox{(direct sum)};\\[1mm]
        (b) &{\pi_{F}}|_G ~\mbox{is an isomorphism from}~ G~ \mbox{onto}~ F;\\[2mm]
\end{array}~~~~~~~~~~~~~~~\quad\quad\quad\quad\]
and that, if they are verified, the projection $E\to G$ associated to the  direct
sum is $({\pi_{F}}|_G)^{-1}\circ \pi_F$. To each  $\alpha \in  L(F;F^{\bot})$ we associate
 its graphic $G=\{x+\alpha(x)\}_{x\in F}$, which is a closed vector subspace of $E$
 verifying the conditions above. Inversely, for each closed vector subspace $G\subset E$ 
verifying the conditions above, there exists one and only one $\alpha\in L(F;F^{\bot})$ whose  graphic is $G$, namely $\alpha=\pi_{F^{\bot}}\circ ({\pi_{F}}|_G)^{-1}$.

        We will use the preceding well-known considerations in the proof of the 
following result:

\subsection{} Let $E$ be a real Hilbert space and let $F\subset E$ be a closed vector 
subspace. Let $\mathcal{U}_F\subset G(E)$ be the set of the orthogonal projections $\xi \in
 G(E)$ such  that $E = F^{\bot}\oplus \xi(E)$. Then $\mathcal{U}_F$ is an open subset in
 $G(E)$, containing $\pi_F$, and 
there exists a diffeomorphism  $\psi_F:\mathcal{U}_F\to L(F;F^{\bot})$, defined by 
$\psi_F(\xi)=\pi_{F^{\bot}}\circ({\pi_{F}}|_{\xi(E)})^{-1}$, that verifies $\psi_F(\pi_F)=0$.

\begin{proof} The considerations before the statement show that $\psi_F$ is a bijective 
map from $\mathcal{U}_F$ onto $L(F;F^{\bot})$, whose inverse $\psi_F^{-1}:L(F;F^{\bot})\to
 \mathcal{U}_F$ associates to each  $\alpha$ the orthogonal projection onto the closed
 vector subspace $\{x+\alpha(x)\}_{x\in F}$.
All we have to show is that $\mathcal{U}_F$ is open in $G(E)$ and that both $\psi_F$ and
 $\psi_F^{-1}$ are smooth maps. For that, we consider the morphism from the tautological
 vector bundle $(\xi(E))_{\xi\in G(E)}$ into the constant vector bundle $F_{G(E)}$
 whose value at $\xi\in G(E)$ is ${\pi_{F}}|_{\xi(E)}:\xi(E)\to F$; the fact that 
$\xi \in \mathcal{U}_F$ if and  only if the ''fibre'' of the morphism at $\xi$ is an
 isomorphism implies that $ \mathcal{U}_F$ is open in $G(E)$; taking the restrictions of 
the vector bundles to $ \mathcal{U}_F$, the  fact that the inverse of a (smooth) isomorphism
 is smooth implies that the map $ \mathcal{U}_F\to L(F;E)$, 
$\xi \to  ({\pi_{F}}|_{\xi(E)})^{-1}$ is smooth, hence $\psi_F:  \mathcal{U}_F
\to L(F;F^{\bot})$ is also smooth. Now, we have an injective morphism from the constant
 vector bundle $F_{L(F,F^{\bot})}$ into the constant vector bundle $E_{L(F,F^{\bot})}$,
 whose fibre at $\alpha \in L(F,F^{\bot})$ is the linear map $F\to E$,  
$x \to x + \alpha(x)$, hence the image of this morphism is a vector bundle with basis
 $L(F;F^{\bot})$ and this implies that the map $\psi_F^{-1}:L(F,F^{\bot}) \to L(E;E)$ is
 smooth.
\end{proof}

   As a corollary, we have:
\subsection{}  If $E$ is a real Hilbert space, then $G(E)$ is a manifold in $L(E;E)$.
 If $E$ is $N$-dimensional and $F\subset E$ is $n$-dimensional, then the dimension 
of $G(E)$ at $\pi_F$ is $n(N-n)$.

\subsection{} Let $E$ be a real Hilbert space, $F\subset E$ be a closed vector subspace
 and $\psi_F:\mathcal{U}_F\to  L(F, F^{\bot})$ be the diffeomorphism defined in 3.4. 
For each $\xi\in \mathcal{U}_F$ and $\eta\in T_{\xi}(G(E))$, we have
\[ D\psi_F(\xi)(\eta) ~=~\eta\circ({\pi_{F}}|_{\xi(E)})^{-1}-({\pi_{F}}|_{\xi(E)})^{-1}
\circ \pi_F\circ \eta\circ({\pi_{F}}|_{\xi(E)})^{-1}. \]
In particular, $D\psi_F(\xi)(\eta)=\eta|_F$.
\begin{proof}
 Let $\phi_F:\mathcal{U}_F\to  L(F;E)$ be the smooth map defined by 
$\phi_F(\xi)=({\pi_{F}}|_{\xi(E)})^{-1}$ (see the proof of 3.4).
 Let $w\in F$ arbitrary. Differentiating the identity
$\pi_F(\phi_F(\xi)(w)) = w$, we obtain 
\[\pi_F(D\phi_F(\xi)(\eta)(w)) = 0,\]
hence $D\phi_F(\xi)(\eta)(w)\in F^{\bot}$. On the other hand, we have a smooth section of
 the tautological vector bundle $(\xi(E))_{\xi\in G(E)}$ associating to each  $\xi$,
 $\phi_F(\xi)(w)$;  its covariant derivative with respect to the metric connection, which,
 by (2.5)  and 3.3, is equal to
\[D\phi_F(\xi)(\eta)(w)- \eta(\phi_F(\xi)(w)),\]
must hence belong to $\xi(E)$. We can now conclude that $D\phi_F(\xi)(\eta)(w)$ is the
 projection of $\eta(\phi_F(\xi)(w))$ onto $F^{\bot}$ associated to the direct sum  
$E = F^{\bot}\oplus \xi(E)$. 
The fact that $\psi_F(\xi)(w)=\pi_{F^{\bot}}(\phi_F(\xi)(w))$ shows that 
$D\psi_F(\xi)(\eta)(w)=\pi_{F^{\bot}}(D\phi_F(\xi)(w))$, 
 hence $D\psi_F(\xi)(\eta)(w)$ is also the projection of $\eta(\phi_F(\xi)(w))$
onto $F^{\bot}$ associated to the direct sum $E = F^{\bot}\oplus \xi(E)$ and, by the
 considerations  made before 3.4, this projection is equal to
\[\eta\left(({\pi_{F}}|_{\xi(E)})^{-1}(w)\right)-({\pi_{F}}|_{\xi(E)})^{-1}
\left(\pi_F(\eta(({\pi_{F}}|_{\xi(E)})^{-1}(w)))\right).\] 
To show that $D\psi_F(\pi_F)(\eta)=\eta|_F$ it will be enough to know that each 
$\eta \in T_{\pi_F} (G(E))$ maps $F$ into $F^{\bot}$.  To see this, we differentiate the identity $\xi\circ\xi=\xi$ and obtain $\eta\circ \pi_F+ \pi_F\circ \eta=\eta$, hence 
$\eta\circ \pi_F=\eta -\pi_F\circ \eta=\pi_{F^{\bot}}\circ \eta$ and the proof is
complete.
\end{proof}

       We present now several equivalent characterizations of the tangent vector 
spaces to G(E).

\subsection{} Let $E$ be a real Hilbert space and let $F\subset E$ be a closed vector 
subspace. The tangent vector space $T_{\pi_F}(G(E))$ is then contained in the vector 
space $L_{sa}(E;E)$ of self adjoint maps and, for each $\eta\in L_{sa}(E;E)$, the following
 conditions are equivalent:
\[\begin{array}{cl}
   \mbox{(a)}& \eta\in T_{\pi_F}(G(E));\\[1mm]
   \mbox{(b)}& \eta(F)\subset F^{\bot}~~\mbox{and}~~ \eta(F^{\bot})\subset F;\\[1mm]
   \mbox{(c)}& \eta\circ \pi_F + \pi_F\circ \eta = \eta;\\[1mm]
   \mbox{(d)}& \eta\circ \pi_F = (Id- \pi_F)\circ \eta ;\\[1mm]
   \mbox{(e)}& \eta\circ (Id-\pi_F) = \pi_F\circ \eta;\\[1mm]
   \mbox{(f)}& \eta\circ (2\pi_F-Id)= -(2\pi_F -Id) \circ \eta.
\end{array}\hspace{5cm}\]
  
\begin{proof} 
The fact that each $T_{\pi_F}(G(E))$ is contained in $L_{sa}(E;E)$ is a consequence 
of the fact that $G(E)\subset L_{sa}(E;E)$. The equivalence between the four last 
conditions is trivial. Assuming (a), we obtain (c) simply by differentiating
the identity $\xi\circ \xi=\xi$ at $\pi_F$ in the direction of $\eta$. It is readily seen
 that condition (d) implies that $\eta(F)\subset F^{\bot}$ and that condition (e) implies
 that  $\eta(F^{\bot})\subset F$ ($Id-\pi_F=\pi_{F^{\bot}}$). Let us prove now that
 condition (b) implies condition (a). The fact that $\psi_F$ is a diffeomorphism from the
 open set $\mathcal{U}_F$ in $G(E)$  onto $L(F;F^{\bot})$ implies that
 $D\psi_F(\pi_F):T_{\pi_F}(G(E))\to L(E;F^{\bot})$ is an isomorphism. We
can hence take $\eta' \in T_{\pi_F}(G(E))$ such that 
\[\eta|_F=D\psi_F(\pi_F)(\eta')=\eta'|_F.\]
Then $\eta'$ is self-adjoint and verifies condition (b), hence 
$\eta'|_{F^{\bot}}:F^{\bot}\to F$ is the adjoint map to $\eta'|_{F}:F\to F^{\bot}$ and 
$ \eta|_{F^{\bot}}:F^{\bot}\to F$ is the adjoint map to $\eta|_{F}:F\to F^{\bot}$. We 
deduce now that    $\eta'|_{F^{\bot}}=\eta|_{F^{\bot}}$, hence $\eta=\eta'$ and the proof 
is complete.
\end{proof}

\noindent
{\em Remark}. To feel what is happening, assume that $E$ is finite dimensional and 
take an orthonormal basis $x_1,\ldots, x_N$ of $E$, whose first $n$ vectors constitute 
a basis for $F$. Then the matrices of $\pi_F$, $id-\pi_F$  and $2\pi_F-Id$ are respectively
\[\left[\begin{array}{cc}
Id&0\\0&0\end{array}\right]
\quad\quad \left[\begin{array}{cc}
0&0\\0&Id\end{array}\right]
\quad\quad \left[\begin{array}{cc}
Id&0\\0&-Id\end{array}\right]
\hspace{5cm}\]
and condition (b) says that the elements of $T_{\pi_F}(G(E))$ are the linear maps
whose matrix has the form
\[ \left[\begin{array}{cc} 0& A^*\\A&0\end{array}\right].\hspace{10cm}\]

\section{The Differential Geometry of Grassman manifolds}

\subsection{} Let $E$ be a real Hilbert space and let $F\subset E$ be a closed vector 
subspace. For each $\eta\in L_{sa}(E;E)$ the following conditions are then equivalent:
\[\begin{array}{cl}
   \mbox{(a)}& \eta(F)\subset F~~\mbox{and}~~ \eta(F^{\bot})\subset F^{\bot};\\[1mm]
   \mbox{(b)}& \eta\circ \pi_F = \pi_F\circ \eta ;\\[1mm]
   \mbox{(c)}& \eta\circ (Id-\pi_F) = (Id-\pi_F)\circ \eta.
\end{array}\hspace{5cm}\]

 We will denote by $T_{\pi_F}(G(E))^{\bot}$ the set of self-adjoint linear maps
$\eta\in  L_{sa}(E;E)$ verifying the preceding conditions.
\begin{proof} The fact that (b) and (c) are equivalent is trivial. It is readily seen that
 (b) implies $\eta(F)\subset F$ and that (c) implies $\eta(F^{\bot})\subset F^{\bot}$. Assuming (a), one sees that $\eta\circ \pi_F(x) = \eta(x) = \pi_F\circ \eta(x)$ for $x\in F$
 and $\eta \circ \pi_F(x)=0= \pi_F\circ \eta(x)$ for $x\in F^{\bot}$, hence 
$\eta\circ \pi_F(x)=\pi_F\circ \eta(x)$ for arbitrary $x$ and (b) is proved.
\end{proof}

\subsection{} Let $E$ be a real Hilbert space and let $F\subset E$ be a closed vector
 subspace. Then $L_{sa}(E;E)$ is the direct sum of the closed vector subspaces 
$T_{\pi_F}(G(E))$  and $T_{\pi_F}(G(E))^{\bot}$ and the projections 
$\bar{\pi}_{\pi_F} : L_{sa}(E;E)\to T_{\pi_F}(G(E))$ and 
$\bar{\pi}_{\pi_F}^{\bot} :L_{sa}(E;E)\to T_{\pi_F}(G(E))^{\bot}$
associated to this direct sum are defined by 
\[\begin{array}{l}
\bar{\pi}_{\pi_F}(\eta)=(Id-\pi_F)\circ \eta \circ \pi_F + \pi_F\circ \eta \circ (Id-\pi_F),\\[1mm]
\bar{\pi}_{\pi_F}^{\bot}(\eta)=(Id-\pi_F)\circ \eta \circ(Id- \pi_F) + \pi_F\circ \eta \circ 
\pi_F.
\end{array}\hspace{0cm}\].
\begin{proof}
Conditions (a) of 4.1 and (b) of 3.7 show that the intersection
$T_{\pi_F}(G(E))\cap T_{\pi_F}(G(E))^{\bot}$ is $\{0\}$. It is readily seen that, for each
$\eta\in L_{sa}(E;E)$, $\bar{\pi}_{\pi_F}(\eta)$ applies $F$ into $F^{\bot}$ and $F^{\bot}$
into $F$ and $\bar{\pi}^{\bot}_{\pi_F}(\eta)$ applies $F$ into $F$ and  $F^{\bot}$ into
 $F^{\bot}$, hence $\bar{\pi}_{\pi_F}(\eta)\in T_{\pi_F}(G(E))$ and
 $\bar{\pi}^{\bot}_{\pi_F}(\eta)\in T_{\pi_F}(G(E))^{\bot}$.
All we have to note now is that, for each $\eta$, $\bar{\pi}_{\pi_F}(\eta)+
\bar{\pi}^{\bot}_{\pi_F}(\eta)=\eta$.
\end{proof}

\subsection{} If $E$ is a finite dimensional real Hilbert space and if we consider in 
$L_{sa}(E;E)$ the Hilbert-Schmidt inner product, then, for each vector subspace 
$F\subset E$, the subspaces $T_{\pi_F}(G(E))$ and $T_{\pi_F}(G(E))^{\bot}$ of $L_{sa}(E;E)$
 are mutually orthogonal, hence each one is the orthgonal complement of the other.

\begin{proof}
Assume $\eta \in T_{\pi_F}(G(E))$ and $\eta' \in T_{\pi_F}(G(E))^{\bot}$.  Choose an
 orthonormal basis $x_1, \ldots, x_N$ of $E$ such that the first $n$ vectors constitute
a basis of $F$ and the last $N-n$ vectors constitute a basis of $F^{\bot}$. Conditions 
(b) of 3.7 and (a) of 4.1 assure that, for each $1 \leq k\leq N$, $\langle \eta(x_k), 
\eta'(x_k)\rangle = 0$, hence $\langle \eta, \eta' \rangle =0$ ( cf.\ (2.2)).
\end{proof}
      
 The preceding result explains why we employ the notation $ T_{\pi_F}(G(E))^{\bot}$
and $\bar{\pi}_{\pi_F}$.\\

 If $E$ is a finite or infinite dimensional real Hilbert space we will define 
the \em canonical connection \em in the manifold $G(E)$ as the one that verifies the 
condition that $\theta_{\pi_F}(\eta, \alpha)$ belongs to the kernel $T_{\pi_F}(G(E))^{\bot}$
of the linear map $\bar{\pi}_{\pi_F}:L_{sa}(E;E)\to T_{\pi_F}(G(E))$, for each $\eta$ and 
$\alpha$ in $ T_{\pi_F}(G(E))$. In an analogous way to that used in the case of a metric
 connection, it is easily seen that this  connection is symmetric and also is defined by 
the formula
\begin{equation}
\theta_{\pi_F}(\eta, \alpha)=D\bar{\pi}_{\pi_F}(\alpha)(\eta).
\end{equation}

This is the connection that we will always consider in the Grassman manifold $G(E)$. 
Of, course, in case $E$ is finite dimensional, this connection is the metric connection 
with respect to the Hilbert-Schmidt inner product.

  We can obtain a more explicit formula for the connection on $G(E)$ by calculating the
 derivative in (4.1), using the formula in 4.2, $\bar{\pi}_{\xi}(\eta)=
(Id-\xi)\circ \eta\circ \xi + \xi \circ \eta \circ (Id-\xi)$. This gives
\begin{equation}\begin{array}{lcl}
\theta_{\xi}(\eta, \alpha)&=& -\alpha\circ \eta\circ\xi +(Id-\xi)\circ \eta\circ \alpha
+\alpha\circ \eta \circ (Id-\xi) -\xi\circ\eta\circ \alpha,\\[1mm]
 &=&(Id-2\xi)\circ \eta\circ \alpha
+\alpha\circ \eta \circ (Id-2\xi).
\end{array}
\end{equation}

Let us now obtain, using (2.7), two formulas for the curvature, the first
for the metric connection of the tautological vector bundle, and the second 
for the canonical connection of the Grassman manifold. In the first case, 
we take $\hat{\theta}_{\xi}(w,\eta) = \eta(w)$ for each $w\in E$ and $\eta \in L(E;E)$
 (cf.\ 3.3) obtaining
\begin{equation}
  R_{\xi}(\alpha,\beta,w)=\beta(\alpha(w))-\alpha(\beta(w))                  
\end{equation}
for each $\alpha$ and $\beta$ in $T_{\xi}(G(E))$ and  $w\in \xi(E)$. In the second case,
we take $\hat{\theta}_{\xi}(\eta, \alpha)=(Id-2\xi)\circ \eta\circ \alpha
+\alpha\circ \eta \circ (Id-2\xi)$ and obtain, noting that
$(Id-2\xi) \circ(Id-2\xi) = Id$ and that, by 3.7(f), $Id-2\xi$ commutes with the composite 
of any two elements of $T_{\xi}(G(E))$,
\begin{equation}
  R_{\xi}(\alpha, \beta, \eta) = \eta\circ\alpha\circ\beta - \eta\circ\beta\circ\alpha
+ \beta\circ\alpha\circ \eta- \alpha\circ\beta\circ \eta.
\end{equation}

 Assuming that E is finite dimensional, we obtain, for the sectional curvatures:
\begin{equation}
\mathrm{Riem}_{\xi}(\alpha,\beta) =  \langle R_{\xi}(\alpha,\beta, \alpha),\beta\rangle = 2\langle \alpha \circ \beta, \alpha\circ \beta\rangle
-2 \langle \alpha \circ \beta, \beta \circ \alpha\rangle.
\end{equation}

To prove it, all we have to do is to apply the formulas in (2.3), 
remembering that $\alpha$ and $\beta$ are self-adjoint.  
The fact that $(\alpha\circ\beta)^*=\beta\circ\alpha$ implies that  $\alpha\circ \beta$
 and $ \beta\circ\alpha$  have the same norm and we can hence apply Cauchy-Schwartz to
conclude that $\mathrm{Riem}_{\xi}(\alpha,\beta)\geq 0$ and 
$\mathrm{Riem}_{\xi}(\alpha,\beta)= 0$ if and only if 
$\alpha\circ \beta =\beta\circ\alpha$.
                                                                         
 One can also establish easily the following formula for the Ricci curvature:
\begin{equation}              
\mathrm{Ricci}_{\xi}(\alpha,\beta)   =\frac{N-2}{2}\langle \alpha,\beta\rangle
\end{equation}    
where $N$ is the dimension of $E$.

 Grassman manifolds (or, more precisely, their connected components) are 
sometimes represented as homogeneous spaces of the orthogonal group. The 
following considerations will compare this approach with the one we are 
using.

 Let $E$ be a real Hilbert space and let $O(E) \subset L(E;E)$ be the orthogonal group,
  i.e. the set of the toplinear isomorphisms $\xi:E\to E$ such that $\xi^*=\xi^{-1}$.
It is well known that $ O(E)$ is a manifold (a Lie group) and that, for each 
   $\xi\in O(E)$ and $\alpha\in  L(E;E)$, we have
\begin{equation}
\alpha \in T_{\xi}(O(E))\, \mbox{~if and only if~} \,\alpha^*\circ \xi +\xi^*\circ\alpha =0.
\end{equation} 

   In the case where $E$ is finite dimensional, the Riemann structure in $O(E)$
induced by the Hilbert-Schmidt inner product is readily seen to be bi-invariant. The
 orthogonal projections $\pi_{\xi} :L(E;E) \to T_{\xi}(O(E))$ are defined by
\begin{equation}
\pi_{\xi}(\lambda)=\frac{1}{2}(\lambda -\xi \circ \lambda^*\circ \xi).
\end{equation}

   Even in the case $E$ is infinite dimensional, we define projection maps 
$\pi_{\xi} :L(E;E) \to T_{\xi}(O(E))$ by formula (4.8) and we have an associated symmetric 
connection in $O(E)$ defined by the bilinear maps 
$\theta_{\xi}: T_{\xi}(O(E))\times  T_{\xi}(O(E))\to L(E;E)$,  
\begin{equation}
 \theta_{\xi}(\alpha,\beta) =D\pi_{\xi}(\beta)(\alpha)=-\frac{1}{2}
\left( \beta\circ \alpha^*\circ \xi + \xi \circ \alpha^*\circ \beta\right).
\end{equation}

   Of course, in the finite dimensional case, this will be the metric connection.

   Now assume that $E$ is a finite or infinite dimensional real Hilbert space 
and $H \subset E$ is a fixed closed vector subspace. We can define a smooth map 
 $\Phi:O(E)\to  G(E)$ associating to each $\xi \in O(E)$ the orthogonal projection onto
$\xi(H)$; denoting by $\pi:E \to H$ the orthogonal projection, it is easy to see that 
we have 
\begin{equation}
\Phi(\xi)= \xi \circ \pi\circ \xi^*.
\end{equation}

   Although $\Phi$ is not a totally geodesic map, we can nevertheless state: 

\subsection{} $\Phi:O(E)\to O(E)$ has totally geodesic fibres and, in case $E$ is finite dimensional, is a Riemannian submersion.

\begin{proof} The derivative linear map 
$D\Phi_{\xi}:T_{\xi}(O(E))\to T_{\Phi(\xi)}(O(E))$ is defined by   
\[ D\Phi_{\xi}(\alpha)= \alpha \circ \pi \circ \xi^* + \xi \circ \pi \circ \alpha^*.\]

Given $\alpha\in T_{\xi}(O(E))$ and $\beta \in T_{\Phi(\xi)}(O(G(E))$ arbitrary, we obtain,
 using (4.7), 3.7(c) and (2.3),
\begin{eqnarray*}
\langle D\Phi_{\xi}(\alpha),\beta\rangle &=& \langle \alpha\circ \pi\circ \xi^*, \beta\rangle
+\langle \xi\circ \pi\circ \alpha^*, \beta\rangle\\
&=& \langle \alpha\circ \pi\circ \xi^*, \beta\rangle
+\langle -\xi\circ \pi\circ \xi^*\circ \alpha\circ \xi^*, \beta\rangle\\
&=& \langle \alpha, \beta\circ\xi\circ \pi\rangle
+\langle \alpha, -\xi\circ \pi\circ \xi^*\circ \beta\circ \xi\rangle\\
&=& \langle \alpha, \beta\circ\xi\circ \pi
-\beta\circ \xi+\beta\circ \xi\circ \pi\rangle \quad = \quad \langle \alpha, 2\beta\circ\xi\circ \pi-\beta\circ \xi\rangle,
\end{eqnarray*}
where, using (4.7), we can see that $2\beta\circ\xi\circ \pi-\beta\circ\xi \in T_{\xi}(O(E))$. Hence, the adjoint linear map $D\Phi_{\xi}^*:T_{\Phi(\xi)}(O(E))\to T_{\xi}(O(E))$
is defined by
\[ D\Phi_{\xi}^*(\beta) = 2\beta\circ \xi \circ \pi -\beta\circ \xi.\]
It is not difficult to verify now that
\[ D\Phi_{\xi}(D\Phi_{\xi}^*(\beta))=\beta,\]
which means precisely that $\Phi$ is a Riemannian submersion.

Let $\xi \in O(E)$ and let $O_o(E)$ be the fibre of $\Phi$ over $\Phi(\xi)$. To prove that
 $O_o(E)$ is a totally geodesic submanifold of $O(E)$, all we have to see is that, for each 
$\alpha$ and $\beta$ in $T_{\xi}(O_o(E))$, we have $(\beta, \theta_{\xi}(\alpha, \beta))\in 
T_{(\xi,\alpha)}(T(O_o(E)))$, where $\theta_{\xi}$ is the  connection on $O(E)$. Using the
 formula for $D\Phi_{\xi}$, we see that the fact that $\alpha$ and $\beta$ are in 
$T_{\xi} (O_o(E))$ is equivalent to
\[ \alpha\circ\pi\circ\xi^* +\xi \circ\pi\circ \alpha^* = 0, \quad
\beta\circ\pi\circ\xi^* +\xi \circ\pi\circ \beta^* = 0\] 
and, using the same formula, one concludes easily that, if $(\beta, \lambda)\in
T_{(\xi, \alpha)}(T(O(E)))$, then $(\beta, \lambda)\in
T_{(\xi, \alpha)}(T(O_o(E)))$ if and only if
\[ \alpha\circ\pi\circ\beta^* +\beta\circ\pi\circ\alpha^*+ \lambda\circ \pi\circ \xi^*
+\xi\circ\pi\circ\lambda^*=0.\]
Now, using formula (4.9) for $\theta_{\xi}(\alpha,\beta)$ and the characterization of
 $T_{\xi}(O(E))$given in 4.7, we obtain
\begin{eqnarray*}
\lefteqn{\alpha\circ\pi\circ\beta^* +\beta\circ\pi\circ\alpha^* +
\theta_{\xi}(\alpha,\beta)\circ\pi\circ\xi^* +\xi\circ\pi\circ\theta_{\xi}(\alpha,\beta)^*}
\\[1mm]
&=& \alpha\circ\pi\circ\beta^* +\beta\circ\pi\circ\alpha^* \\
&&-\frac{1}{2}\left(\beta\circ \alpha^*\circ\xi +\xi\circ \alpha^*\circ \alpha^*\circ \beta
\right)\circ \pi \circ \xi^*- \frac{1}{2}\xi\circ \pi\left(\xi^*\circ \alpha\circ \beta^* + \beta^*\circ\alpha\circ\xi^*\right)\\
&=& \alpha\circ \pi \circ \beta^*+ \beta \circ \pi\circ \alpha^* 
+ \frac{1}{2} \beta\circ\xi^*\circ \alpha\circ \pi\circ\xi^*
+  \frac{1}{2} \alpha\circ\xi^*\circ \beta\circ \pi\circ\xi^*\\
&&+\frac{1}{2}\xi\circ\pi\circ\alpha^*\circ \xi\circ\beta^*
+ \frac{1}{2}\xi\circ\pi\circ\beta^*\circ \xi\circ\alpha^*\\
&=&  \alpha\circ \pi \circ \beta^*+ \beta \circ \pi\circ \alpha^*
 -\frac{1}{2}\beta\circ \pi\circ \alpha^* -\frac{1}{2}\alpha\circ \pi\circ \beta^*
 -\frac{1}{2}\alpha\circ \pi\circ \beta^*-\frac{1}{2}\beta\circ \pi\circ \alpha^*\\
&=& 0,
\end{eqnarray*}
 and the proof is complete.
\end{proof}

We are going now to present a formula for the geodesics in $G(E)$ with arbitrary initial
 conditions. Let $E$ be a real Hilbert space.
\subsection{} For each $\xi\in G(E)$ and $\eta\in T_{\xi}(G(E))$, there exists a smooth map
$f:\mathbb{R}\to G(E)$ defined by
\[ f(t)=\frac{1}{2}\La{(} Id + (2\xi-Id) \circ \cos(2t\eta) + \sin( 2t\eta)\La{)}\]
and $f$ is a geodesic of $G(E)$ that verifies $f(0)=\xi$ and $f'(0)=\eta$.

\begin{proof}
 We note first that, from 3.7(f), we conclude that $(2\xi-Id)$ commutes 
with $\cos(2t\eta)$ and anti-commutes with $\sin(2t\eta)$. It is now trivial that $f(t)$ 
is self-adjoint and, noting that $(2\xi-Id)\circ(2\xi-Id) = Id$ and
\[\cos(2t\eta)\circ \cos(2t\eta) + \sin(2t\eta)\circ\sin(2t\eta) = Id,\]
we obtain
\begin{eqnarray*} 
f(t)\circ f(t) &=& \frac{1}{4}\LA{(}Id + (2\xi-Id)\circ \cos(2t\eta)+\sin(2t\eta)
\quad+(2\xi-Id)\circ \cos(2t\eta)\quad\quad\quad\\
&&\quad+ (2\xi-Id)\circ \cos(2t\eta)\circ (2\xi-Id)\circ \cos(2t\eta)
+(2\xi-Id)\circ \cos(2t\eta)\circ \sin(2t\eta)\\
&&\quad+\sin(2t\eta)+ \sin(2t\eta)\circ (2\xi-Id)\circ \cos(2t\eta)+\sin(2t\eta)\circ \sin(2t\eta)\LA{)}\\[1mm]
&=&  \frac{1}{4}\LA{(}Id + (2\xi-Id)\circ \cos(2t\eta) + \sin(2t\eta)+  (2\xi-Id)\circ \cos(2t\eta)\\
&&\quad+ \cos(2t\eta)\circ \cos(2t\eta) +\sin(2t\eta) +\sin(2t\eta)\circ \sin(2t\eta)
\LA{)} \quad=
\quad f(t),
\end{eqnarray*}
whence we conclude that $f(t) \in G(E)$. Next we see that
\begin{eqnarray*}
f'(t) &=& \frac{1}{2}\La{(} -2(2\xi-Id)\circ \sin(2t\eta)\circ \eta 
+ 2 \cos(2t\eta)\circ\eta\La{)}\\
&=&\La{(} \cos(2t\eta)-(2\xi-Id)\circ \sin(2t\eta)\La{)}\circ \eta, 
\end{eqnarray*}
in particular $f'(0)=\eta$. Next we obtain 
\[f''(t)=\La{(} -2\sin(2t\eta)-2(2\xi-Id)\circ \cos(2t\eta)\La{)}\circ \eta^2.\]
On the other side, remembering 3.7(f), we have
\begin{eqnarray*}
f'(t)\circ f'(t) &=& \cos(2t\eta)^2\circ \eta^2 - \cos(2t\eta)\circ \eta \circ (2\xi-Id)\circ \sin(2t\eta)\circ \eta\\
&& -(2\xi-Id)\circ \sin(2t\eta)\circ \eta\circ \cos(2t\eta)\circ\eta\\
&& + (2\xi-Id)\circ \sin(2t\eta)\circ \eta\circ(2\xi-Id)\circ \sin(2t\eta)\circ\eta\\[1mm]
&=&  \cos(2t\eta)^2\circ \eta^2 +(2\xi-Id)\circ \cos(2t\eta)\circ \sin(2t\eta)\circ\eta^2\\
&&-(2\xi-Id)\circ \sin(2t\eta)\circ \cos(2t\eta)\circ\eta^2\\
&&+(2\xi-Id)\circ (2\xi-Id)\circ \sin(2t\eta)^2\circ\eta^2\\[1mm]
&=& \La{(}  \cos(2t\eta)^2+ \sin(2t\eta)^2\La{)}\circ\eta^2 \quad=\quad \eta^2,
\end{eqnarray*}
and, using (4.2), we have now
\begin{eqnarray*}
\lefteqn{\theta_{f(t)}(f'(t), f'(t))}\\
&=& (Id-2f(t))\circ f'(t)\circ f'(t) +  f'(t)\circ f'(t)\circ (Id-2f(t))\\
&=& \La{(} -(2\xi-Id)\circ\cos(2t\eta)- \sin(2t\eta)\La{)}\circ\eta^2\\
&& + \eta^2\circ \La{(} -(2\xi-Id)\circ\cos(2t\eta)- \sin(2t\eta)\La{)}\quad=\quad f''(t),
\end{eqnarray*}
whence we conclude that $f$ is indeed a geodesic.
\end{proof}

\subsection{} Let $E$ be a real Hilbert space. $G(E)$ is then a symmetric space and, for
 each $\pi\in G(E)$, the symmetry $Sym:G(E)\to G(E)$ with respect to $\pi$ is defined by
\[Sym(\xi)=(Id -2\pi)\circ \xi \circ (Id -2\pi).\]

\begin{proof}
It is trivial that $Sym(\xi)$ is a self-adjoint map and the fact that 
$(Id-2\pi)\circ (Id-2\pi)=Id$ shows that $Sym(\xi)\circ Sym(\xi)=Sym(\xi)$, hence
$Sym(\xi)\in G(E)$. It is trivial to see that $Sym(\pi)=\pi$ and that
$Sym(Sym(\xi))=\xi$. We have
\[ DSym_{\xi}(\alpha)= (Id-2\pi)\circ \alpha\circ (Id-2\pi),\]
hence, remembering (4.2),
\begin{eqnarray*}
\lefteqn{\nabla D Sym_{\xi}(\alpha, \beta) = (Id-2\pi)\circ \theta_{\xi}(\beta, \alpha)\circ (Id-2\pi)- \theta_{Sym(\xi)}(DSym_{\xi}(\beta), DSym_{\xi}(\alpha))}\\[1mm]
&=& (Id-2\pi)\circ (Id-2\xi)\circ \beta\circ \alpha\circ (Id-2\pi)
+ (Id-2\pi)\circ\alpha\circ \beta\circ(Id-2\xi)\circ (Id-2\pi)\\
&&-\La{(} Id -2(Id-2\pi)\circ \xi\circ(Id-2\pi)\La{)}\circ (Id-2\pi)\circ \beta
\circ(Id-2\pi)\circ(Id-2\pi)\circ \alpha\circ(Id-2\pi)\\
&&-(Id-2\pi)\circ \alpha\circ(Id-2\pi)\circ(Id-2\pi)\circ \beta\circ (Id-2\pi)\circ
\La{(} Id -2(Id-2\pi)\circ \xi\circ(Id-2\pi)\La{)}\\[1mm]
&=&(Id-2\pi)\circ \beta\circ \alpha\circ(Id-2\pi)- 2(Id-2\pi)\circ\xi\circ \beta\circ \alpha\circ(Id-2\pi)\\
&&+(Id-2\pi)\circ \alpha\circ\beta\circ (Id-2\pi)- 2(Id-2\pi)\circ\alpha\circ \beta\circ \xi\circ(Id-2\pi)\\
&&-(Id-2\pi)\circ \beta\circ \alpha\circ(Id-2\pi)+ 2(Id-2\pi)\circ\xi\circ \beta\circ \alpha\circ(Id-2\pi)\\
&&-(Id-2\pi)\circ \alpha\circ\beta\circ (Id-2\pi)+ 2(Id-2\pi)\circ\alpha\circ \beta\circ \xi\circ(Id-2\pi)\quad=\quad 0,
\end{eqnarray*}
that is to say, $Sym$ is a totally geodesic diffeomorphism. Now, if 
$f:\mathbb{R}\to  G(E)$ is a geodesic with $f(0) = \pi$ and $f'(0) = \eta$, we have
\[ f(t)=\frac{1}{2}\left( Id + (2\pi-Id)\circ \cos(2t\eta)+ \sin(2t\eta)\right),\]
hence
\begin{eqnarray*}
Sym(f(t)) &=& \frac{1}{2}(Id-2\pi)\circ\La{(} Id + (2\pi-Id)\circ \cos(2t\eta)+
 \sin(2t\eta)\La{)}\circ (Id -2\pi)\\
&=& \frac{1}{2}\La{(} Id + (Id-2\pi)\circ \cos(2t\eta)- \sin(2t\eta)\La{)}\quad =\quad f(-t),
\end{eqnarray*}
and the proof is complete.
\end{proof}
\section{The complex Grassman manifolds}

Assume that $E$ is a complex Hilbert space, whose inner product will always 
be denoted by $\langle,\rangle_{\mathbb{C}}$. Then $E$ is also a real Hilbert space,  with
 the inner product
\begin{equation}
\langle x,y\rangle={\mathrm{Re}}\langle x, y\rangle_{\mathbb{C}}
\end{equation}
and the following two facts are trivial:

\subsection{} If $F\subset E$ is a complex vector subspace, then the orthogonal projection 
$\pi:E\to F$ is the same when we consider in $E$ either the complex or the real 
inner product.

\subsection{}  If $\xi:E \to E$ is a complex linear map, then the adjoint map 
$\xi^*:E\to E$ is the same when we consider $E$ to be either a complex or a real Hilbert space.

   We will denote by $L(E;E)$ the vector space of all continuous \em real \em  linear 
maps and by $L_{\mathbb{C}}(E;E)$ its vector subspace whose elements are the complex linear 
maps. In the case where $E$ is finite dimensional the Hilbert-Schmidt inner 
product that we will consider in $L(E;E)$ will be the one associated to the 
real structure of $E$ and we will consider in the closed subspace $L_{\mathbb{C}}(E;E)$
 the induced inner product.

   If $E$ is a complex Hilbert space, we will denote by $G_{\mathbb{C}}(E)$ the set of the 
orthogonal 
projections onto closed complex vector subspaces, and we call
$G_{\mathbb{C}}(E)$ the \em complex Grassman manifold \em of $E$. $G(E)$ will denote the
 real Grassman manifold of $E$, i.e.\ the Grassman manifold of $E$, when considered as a 
real Hilbert space.  It is trivial to conclude that
\begin{equation}
G_{\mathbb{C}}(E)=G(E)\cap L_{\mathbb{C}}(E;E).
\end{equation}
                                         
   All that has been said in Section 3 applies \em mutatis mutandis \em to the complex
 Grassman manifolds, but one must be aware that $G_{\mathbb{C}}(E)$ is only a \em real \em 
 manifold within the complex vector space $L_{\mathbb{C}}(E;E)$. The essential reason for 
this is the fact that the map $\xi \to \xi^*$ is not $\mathbb{C}$-linear, but it was 
natural to  anticipate this because, in case $E$ is finite dimensional, 
$G_{\mathbb{C}}(E)$ (like $G(E)$) is  compact (because it is closed and bounded) and it is
 well known that there exists no compact nontrivial complex submanifold of a complex vector
 space.

   For each closed complex vector subspace $F\subset E$, we still have a diffeomorphism 
$\psi_F:\mathcal{U}_{F}\to L_{\mathbb{C}}(F, F^{\bot})$, where 
$\mathcal{U}_F$ is open in $G_{\mathbb{C}}(E)$ and contains $\pi_F$ (cf.\ 3.4), hence:

\subsection{} If $E$ has complex dimension $N$ and $F\subset E$ has complex dimension $n$,
 then the real manifold $G_{\mathbb{C}}(E)$ has dimension $2n(N-n)$ in $\pi_F$.

   The tangent vector space $T_{\pi_F}(G_{\mathbb{C}}(E))$ is contained in the real vector
 space $L_{\mathbb{C}sa}(E;E)$, whose elements are the self-adjoint complex linear maps and,
 for each $\eta\in L_{\mathbb{C}sa}(E;E)$, the fact that $\eta \in T_{\pi_F}
(G_{\mathbb{C}}(E))$  is equivalent to each of the 
conditions (b) to (f) of 3.7; in other words:
\begin{equation}
T_{\pi_F}(G_{\mathbb{C}}(E))=T_{\pi_F}(G(E))\cap L_{\mathbb{C}}(E;E).
\end{equation}

Although $G_{\mathbb{C}}(E)$ is only a real submanifold of $ L_{\mathbb{C}}(E;E)$, it admits
a complex structure:
\subsection{} Let $E$ be a complex Hilbert space. Then the real manifold $G_{\mathbb{C}}(E)$
 admits a complex structure defined by the linear maps
\[ J_{\xi}:T_{\xi}(G_{\mathbb{C}}(E))\to T_{\xi}(G_{\mathbb{C}}(E)),\quad
J_{\xi}(\eta)=i\eta\circ (2\xi-Id).\]

For this structure the real diffeomorphisms $\psi_F:\mathcal{U}_{F}\to L_{\mathbb{C}}
(F, F^{\bot})$ are in fact holomorphic.

\begin{proof} To see that $J_{\xi}$ applies $T_{\xi}(G_{\mathbb{C}}(E))$ into itself we 
use 3.7(f), remembering that  $(2\xi -Id)\circ (2\xi-Id)=Id$  and noting that
\[(i\eta\circ (2\xi-Id))^* =-i(\eta\circ (2\xi-Id))^*=-i(2\xi -Id)\circ \eta = i\eta \circ (2\xi -Id).\]
It is also trivial that $J_{\xi}(J_{\xi}(\eta))=-\eta$. The fact that this almost complex 
structure is indeed a complex one comes from the fact that the real diffeomorphisms
$\psi_F:\mathcal{U}_{F}\to L_{\mathbb{C}}(F, F^{\bot})$ are  holomorphic; this is a simple
 consequence of the formula in 3.6,
\[D\psi_F(\xi)(\eta) = \eta \circ (\pi_F|_{\xi(E)})^{-1}-  (\pi_F|_{\xi(E)})^{-1}\circ\pi_F
\circ \eta\circ (\pi_F|_{\xi(E)})^{-1},\]
the formula $D\psi_F(\xi)(J_{\xi}(\eta))=i D\psi_F(\xi)(\eta)$ being a simple consequence 
of the fact that the restriction of $(2\xi-Id)$ to $\xi(E)$  is the identity.
\end{proof}

    Note that, in case $E$ is finite dimensional, if we choose a complex orthonormal basis 
$x_1,\ldots, x_N$ of $E$ such that $x_1,\ldots, x_n$ is a basis of $\xi(E)$, then, if
$\eta\in T_{\xi}(G_{\mathbb{C}}(E))$ has matrix 
$\scr{\left[\begin{array}{cc}0&A^*\\A&0\end{array}\right]}$, $J_{\xi}(\eta)$ has matrix 
$\scr{\left[\begin{array}{cc}0&-iA^*\\iA&0\end{array}\right]}$.\\

 The considerations in 4.1-4.3 and (4.1)-(4.5) apply \em mutatis mutandis \em to the 
complex Grassman manifolds and we have in particular a \em canonical symmetric 
connection \em in $G_{\mathbb{C}}(E)$ defined also by
\begin{equation}
\theta_{\xi}(\eta,\alpha)=(Id-2\xi)\circ\eta\circ\alpha+ \alpha\circ \eta\circ (Id-2\xi).
\end{equation}
This implies in particular that:

\subsection{} $G_{\mathbb{C}}(E)$ is a totally geodesic submanifold of $G(E)$. We see now
 that:
\subsection{} If $E$ is a complex Hilbert space, then the morphism $J= (J_{\xi})$, from the
vector bundle $T(G_{\mathbb{C}}(E))$ into itself, is parallel.

\begin{proof} From $J_{\xi}(\beta)=i\beta \circ ( 2\xi -Id)$, we obtain
\begin{eqnarray*}
\nabla J_{\xi}(\alpha)(\beta) &=& 2i\beta \circ \alpha + i\theta_{\xi}(\beta, \alpha)
\circ (2\xi-Id)-\theta_{\xi}(i\beta \circ(2\xi-Id),\alpha)\\[1mm]
&=& 2i\beta \circ \alpha + i(Id-2\xi)\circ \beta\circ \alpha\circ (2\xi -Id)+ 
i \alpha\circ \beta \circ (Id-2\xi)\circ (2\xi-Id)\\
&&-i(Id -2\xi)\circ\beta \circ(2\xi-Id)\circ\alpha- i\alpha\circ \beta \circ (2\xi -Id)
\circ (Id -2\xi)\\[1mm]
&=& 2i\beta \circ \alpha - i\beta \circ \alpha -i\alpha \circ \beta - i\beta \circ \alpha
+i\alpha \circ \beta \quad =\quad 0.
\end{eqnarray*}
\end{proof}

In the case where the complex Hilbert space $E$ is finite dimensional, we note that the
 Hilbert-Schmidt inner product in $L(E;E)$ is the real part of the complex inner product
(the one defined by (2.2) with $\langle, \rangle_{\mathbb{C}}$ instead of $\langle,
 \rangle$)  and we see that
\begin{eqnarray*}
\langle J_{\xi}(\alpha),J_{\xi}(\beta)\rangle &=& \langle i\alpha\circ (2\xi -Id),
i\beta\circ (2\xi -Id)\rangle \quad = \quad
\langle \alpha\circ (2\xi -Id),
\beta\circ (2\xi -Id)\rangle\\
& = & \langle \alpha\circ (2\xi -Id)\circ 
(2\xi -Id), \beta\rangle \quad = \quad \langle \alpha, \beta\rangle,
\end{eqnarray*}
hence:
\subsection{} If $E$ is a finite dimensional complex Hilbert space, then 
$G_{\mathbb{C}}(E)$ is a K\"{a}hler manifold.

We end with the remark that (4.7)-(4.10) and 4.4 work equally well in the  
complex case, the usual notation for $O_{\mathbb{C}}(E)$ being $U(E)$ (the unitary group).

\clearpage

\oddsidemargin -0.7cm
\evensidemargin -0.7cm
\thispagestyle{empty}
\mbox{ } \\[-3cm]
\mbox{FIGURE 1. Scanned page 85 of the original paper
published in {\it Res.\ Notes Math.} {\bf131} (1985), 85--102}
\mbox{ } \\[2mm]
\includegraphics[width=18cm]{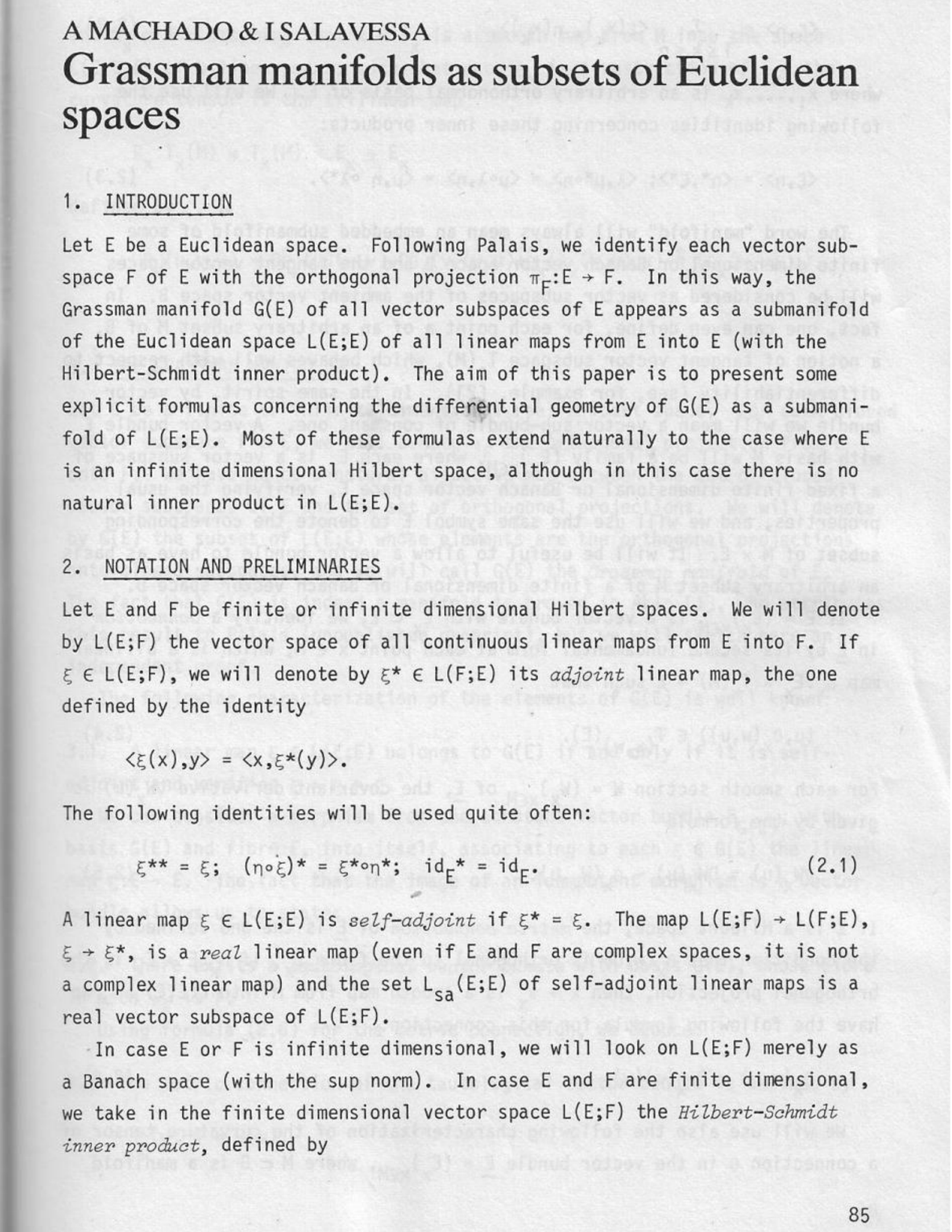}
\vspace*{-7cm}
\clearpage
\thispagestyle{empty}
\mbox{ } \\[-3cm]
\mbox{FIGURE 2. Scanned page 86 of the original paper
published in {\it Res.\ Notes Math.} {\bf131} (1985), 85--102}
\mbox{ } \\[2mm]
\includegraphics[width=18cm]{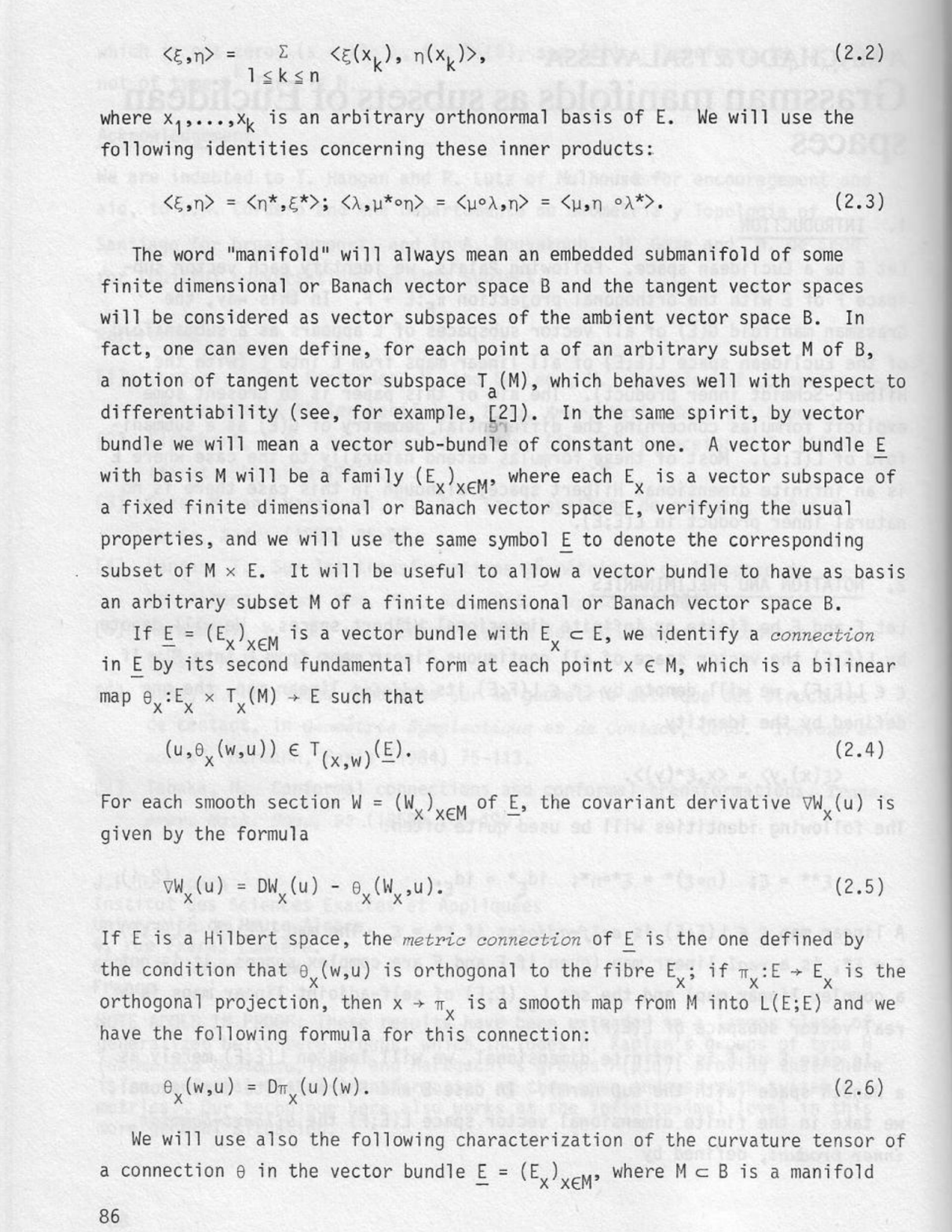}
\vspace*{-7cm}
\clearpage
\thispagestyle{empty}
\mbox{ } \\[-3cm]
\mbox{FIGURE 3. Scanned page 87 of the original paper
published in {\it Res.\ Notes Math.} {\bf131} (1985), 85--102}
\mbox{ } \\[2mm]
\includegraphics[width=18cm]{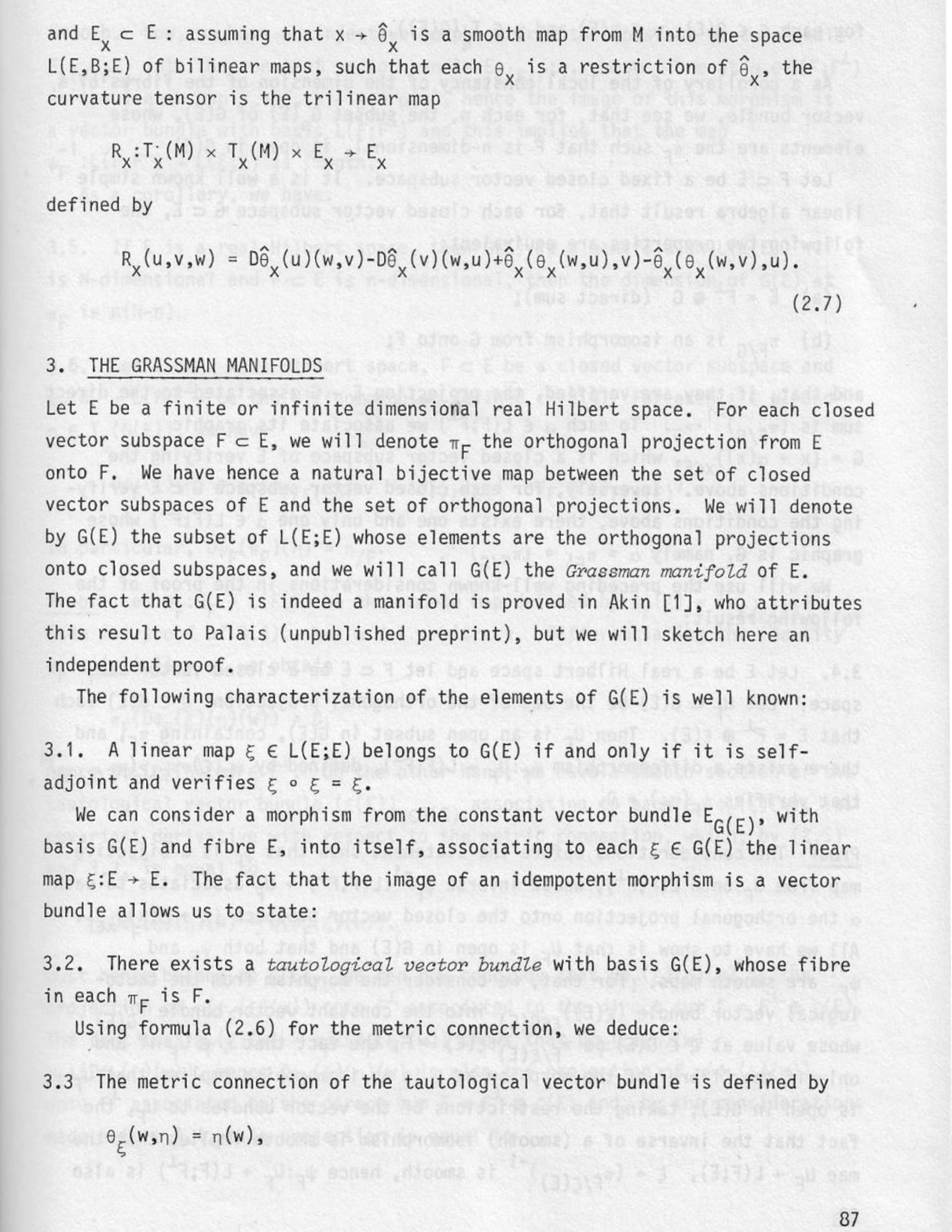}
\vspace*{-7cm}
\clearpage
\thispagestyle{empty}
\mbox{ } \\[-3cm]
\mbox{FIGURE 4. Scanned page 88 of the original paper
published in {\it Res.\ Notes Math.} {\bf131} (1985), 85--102}
\mbox{ } \\[2mm]
\includegraphics[width=18cm]{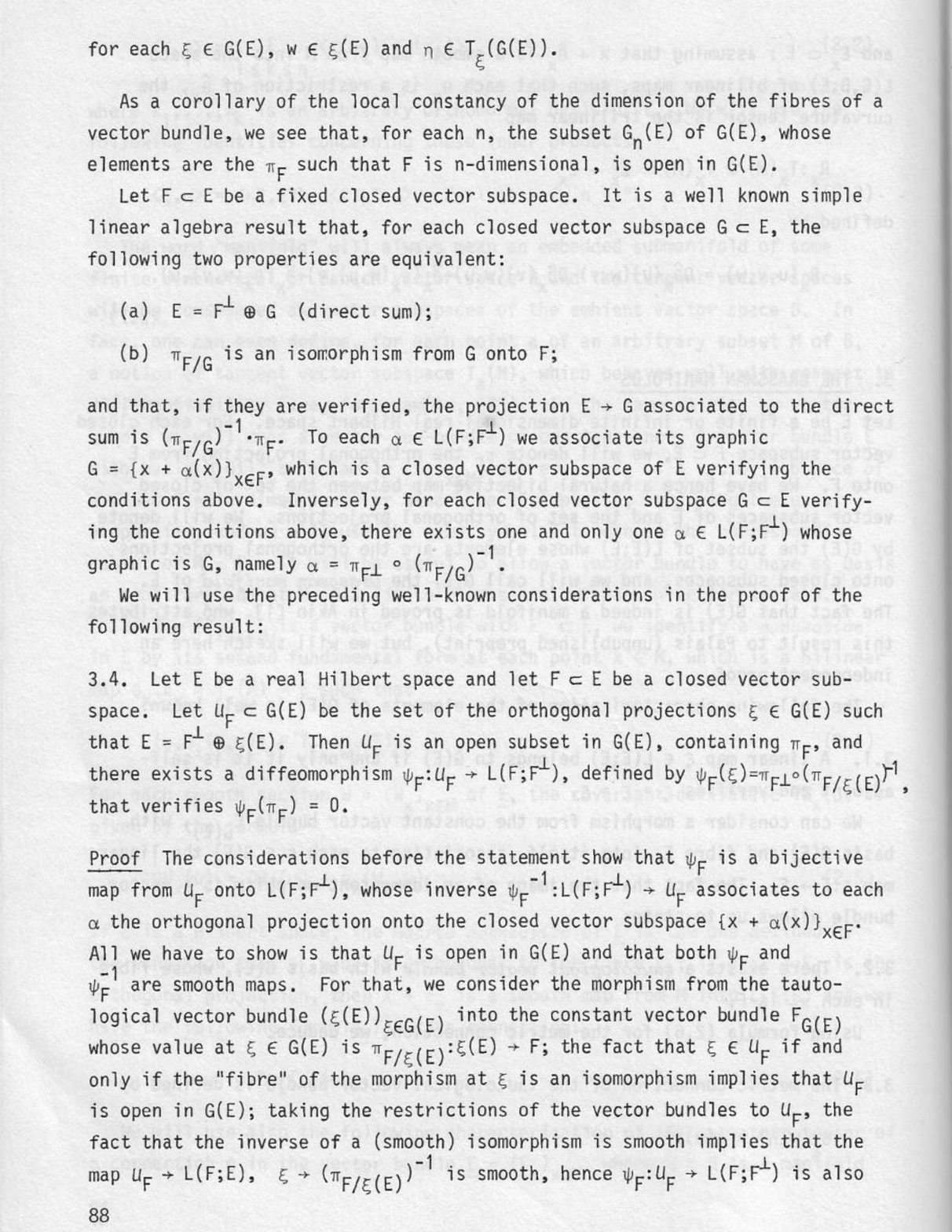}
\vspace*{-7cm}
\clearpage
\thispagestyle{empty}
\mbox{ } \\[-3cm]
\mbox{FIGURE 5. Scanned page 89 of the original paper
published in {\it Res.\ Notes Math.} {\bf131} (1985), 85--102}
\mbox{ } \\[2mm]
\includegraphics[width=18cm]{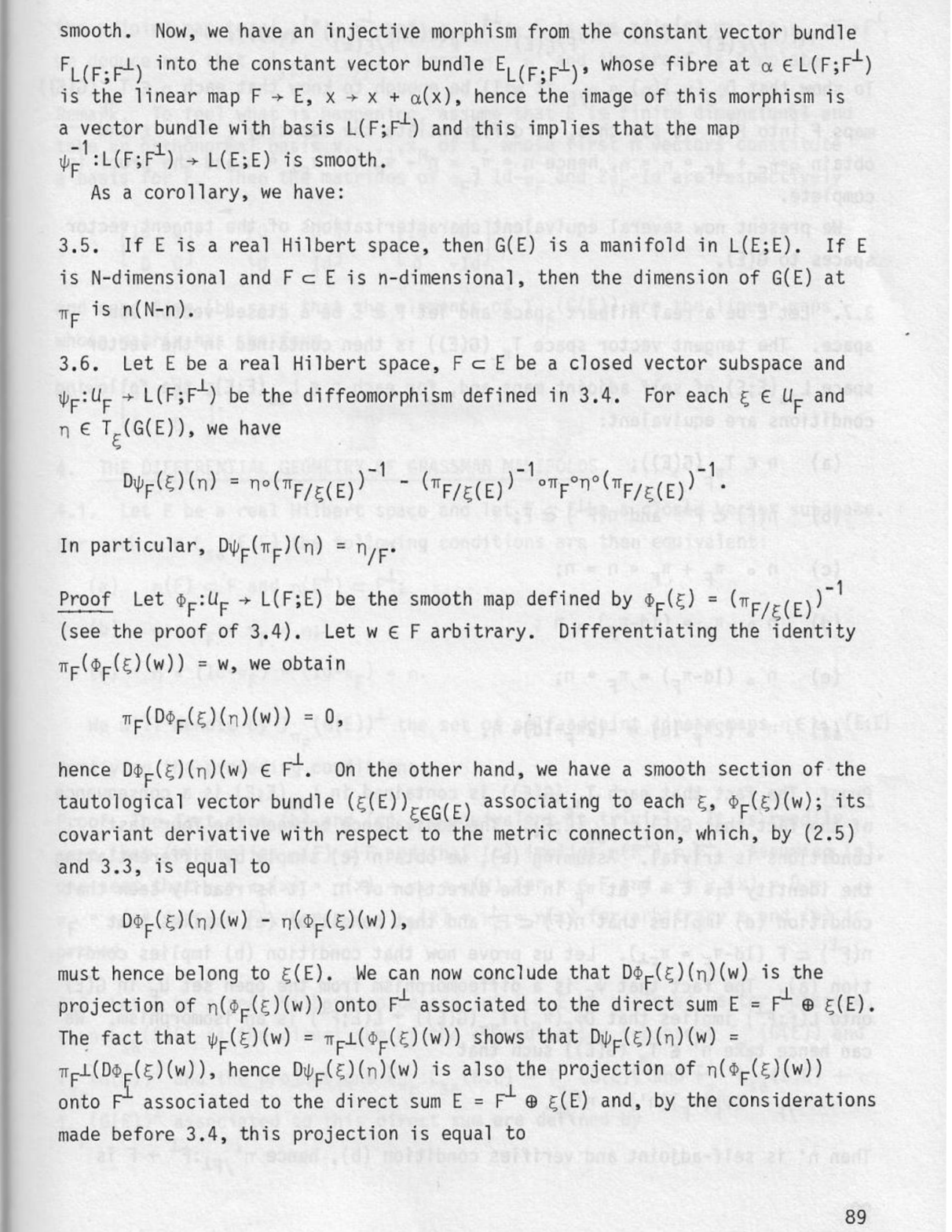}
\vspace*{-7cm}
\clearpage
\thispagestyle{empty}
\mbox{ } \\[-3cm]
\mbox{FIGURE 6. Scanned page 90 of the original paper
published in {\it Res.\ Notes Math.} {\bf131} (1985), 85--102}
\mbox{ } \\[2mm]
\includegraphics[width=18cm]{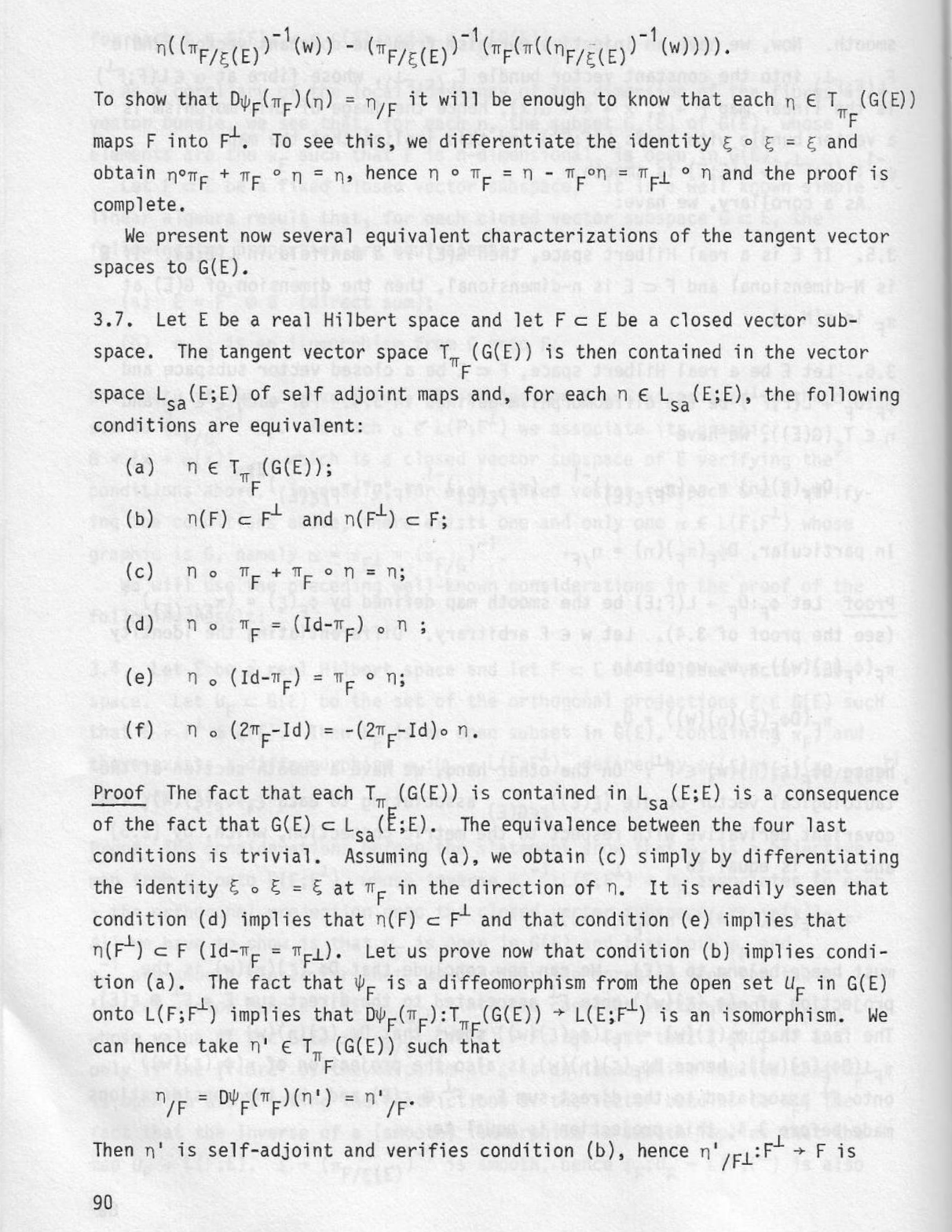}
\vspace*{-7cm}
\clearpage
\thispagestyle{empty}
\mbox{ } \\[-3cm]
\mbox{FIGURE 7. Scanned page 91 of the original paper
published in {\it Res.\ Notes Math.} {\bf131} (1985), 85--102}
\mbox{ } \\[2mm]
\includegraphics[width=18cm]{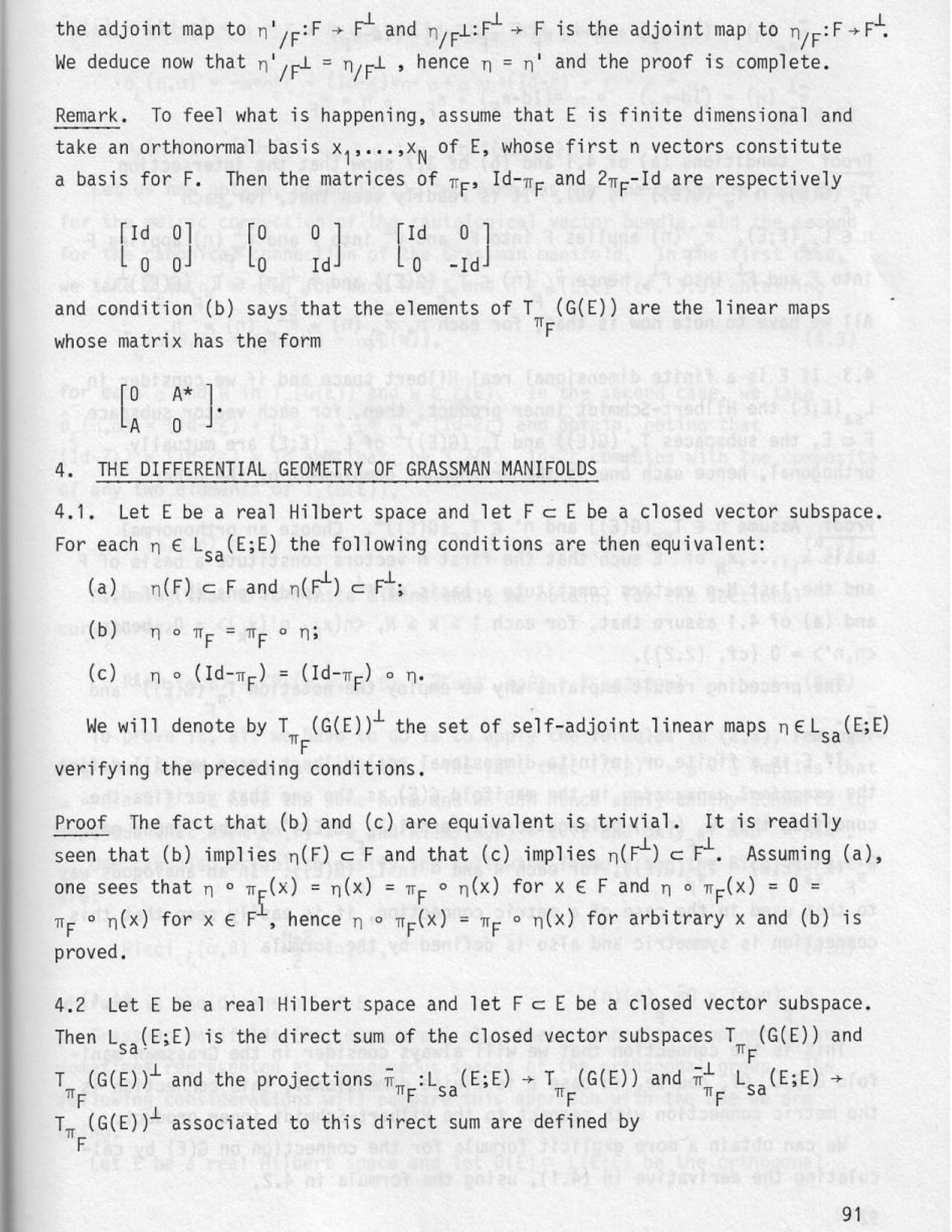}
\vspace*{-7cm}
\clearpage
\thispagestyle{empty}
\mbox{ } \\[-3cm]
\mbox{FIGURE 8. Scanned page 92 of the original paper
published in {\it Res.\ Notes Math.} {\bf131} (1985), 85--102}
\mbox{ } \\[2mm]
\includegraphics[width=18cm]{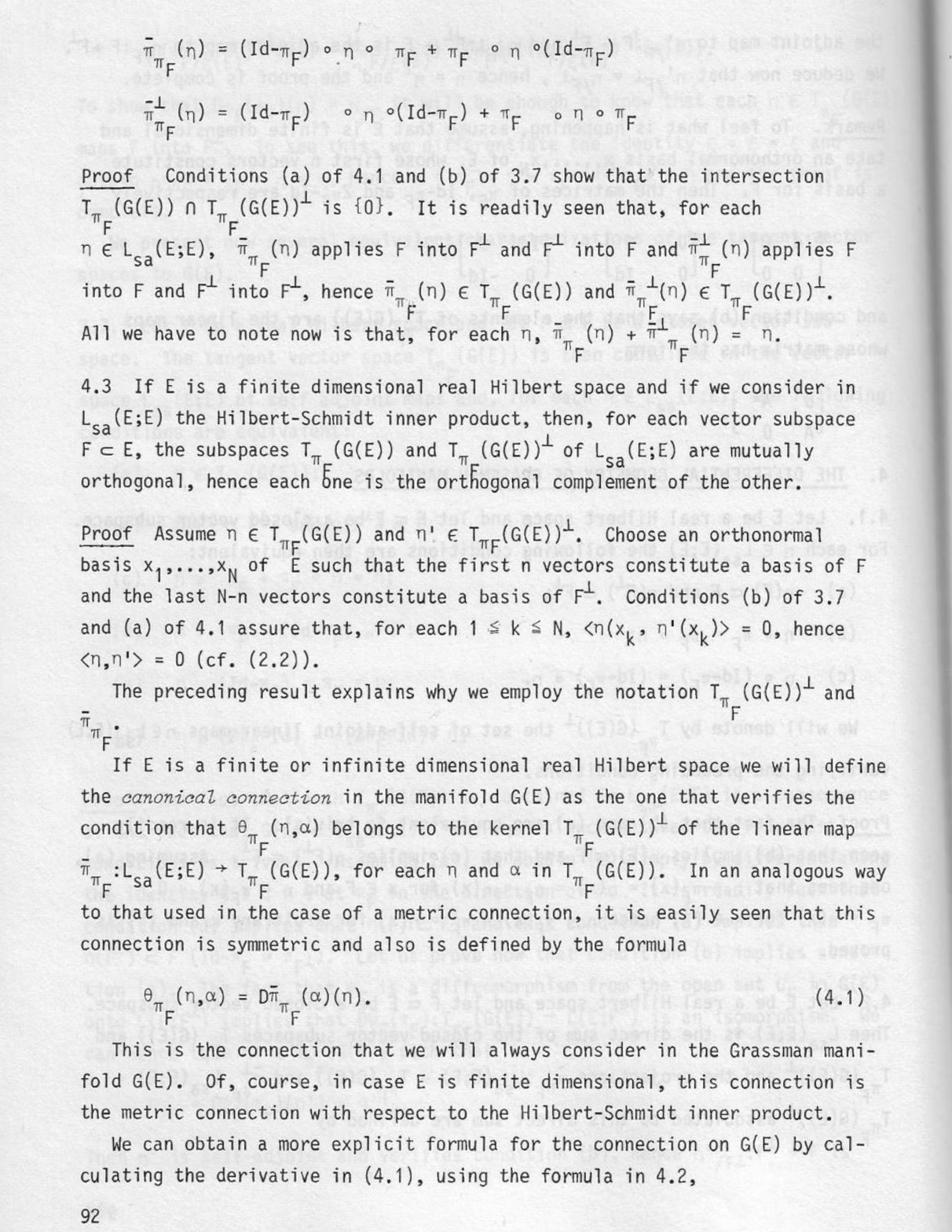}
\vspace*{-7cm}
\clearpage
\thispagestyle{empty}
\mbox{ } \\[-3cm]
\mbox{FIGURE 9. Scanned page 93 of the original paper
published in {\it Res.\ Notes Math.} {\bf131} (1985), 85--102}
\mbox{ } \\[2mm]
\includegraphics[width=18cm]{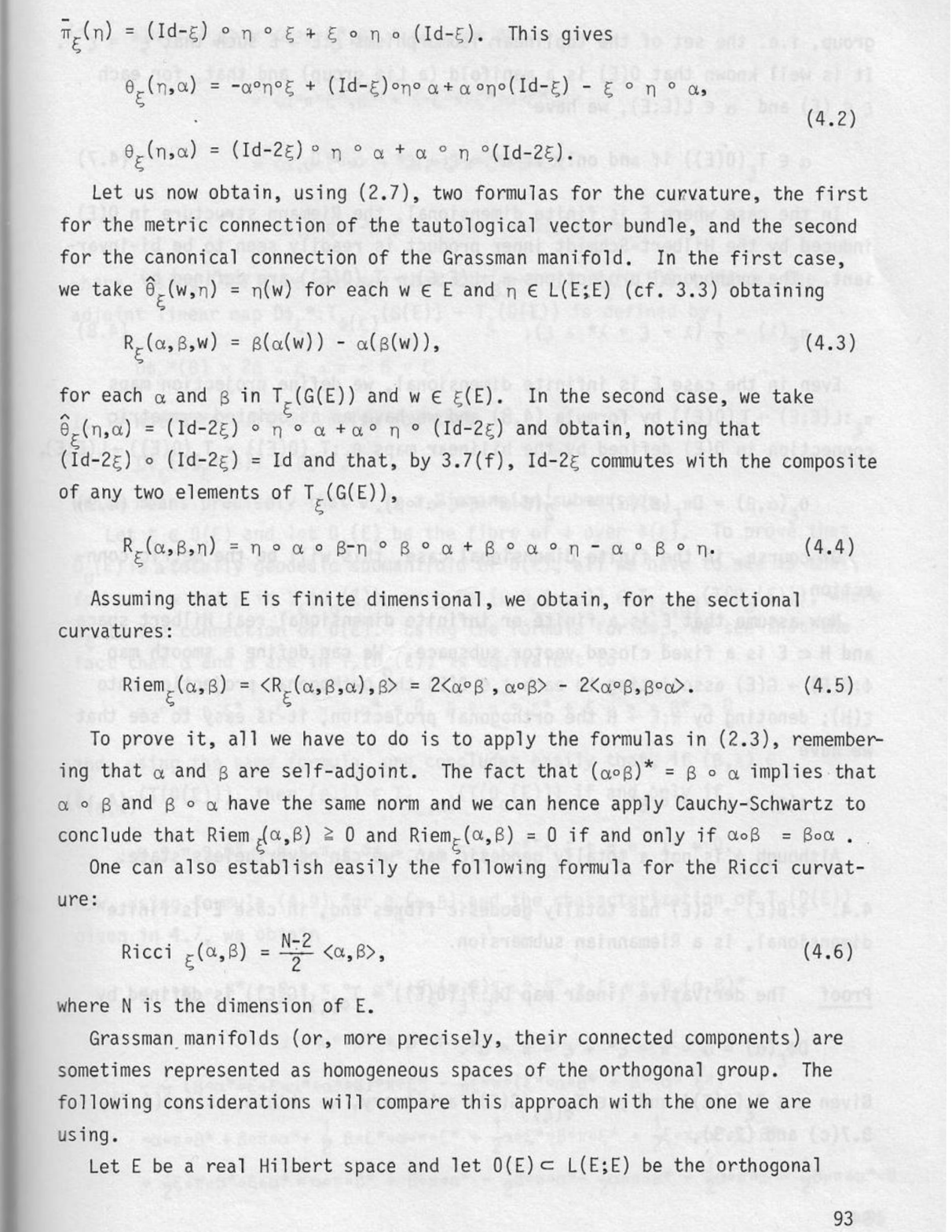}
\vspace*{-7cm}
\clearpage
\thispagestyle{empty}
\mbox{ } \\[-3cm]
\mbox{FIGURE 10. Scanned page 94 of the original paper
published in {\it Res.\ Notes Math.} {\bf131} (1985), 85--102}
\mbox{ } \\[2mm]
\includegraphics[width=18cm]{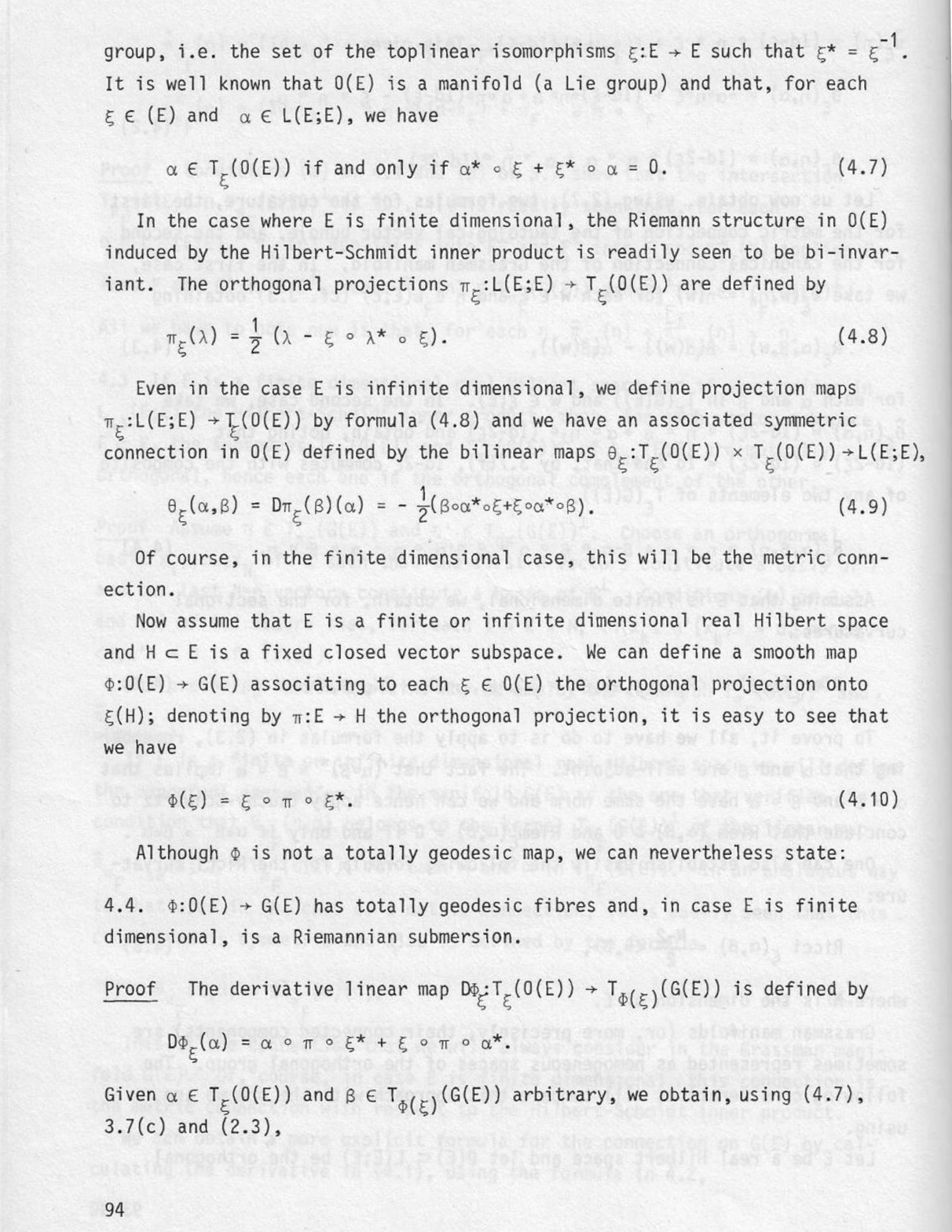}
\vspace*{-7cm}
\clearpage
\thispagestyle{empty}
\mbox{ } \\[-3cm]
\mbox{FIGURE 11. Scanned page 95 of the original paper
published in {\it Res.\ Notes Math.} {\bf131} (1985), 85--102}
\mbox{ } \\[2mm]
\includegraphics[width=18cm]{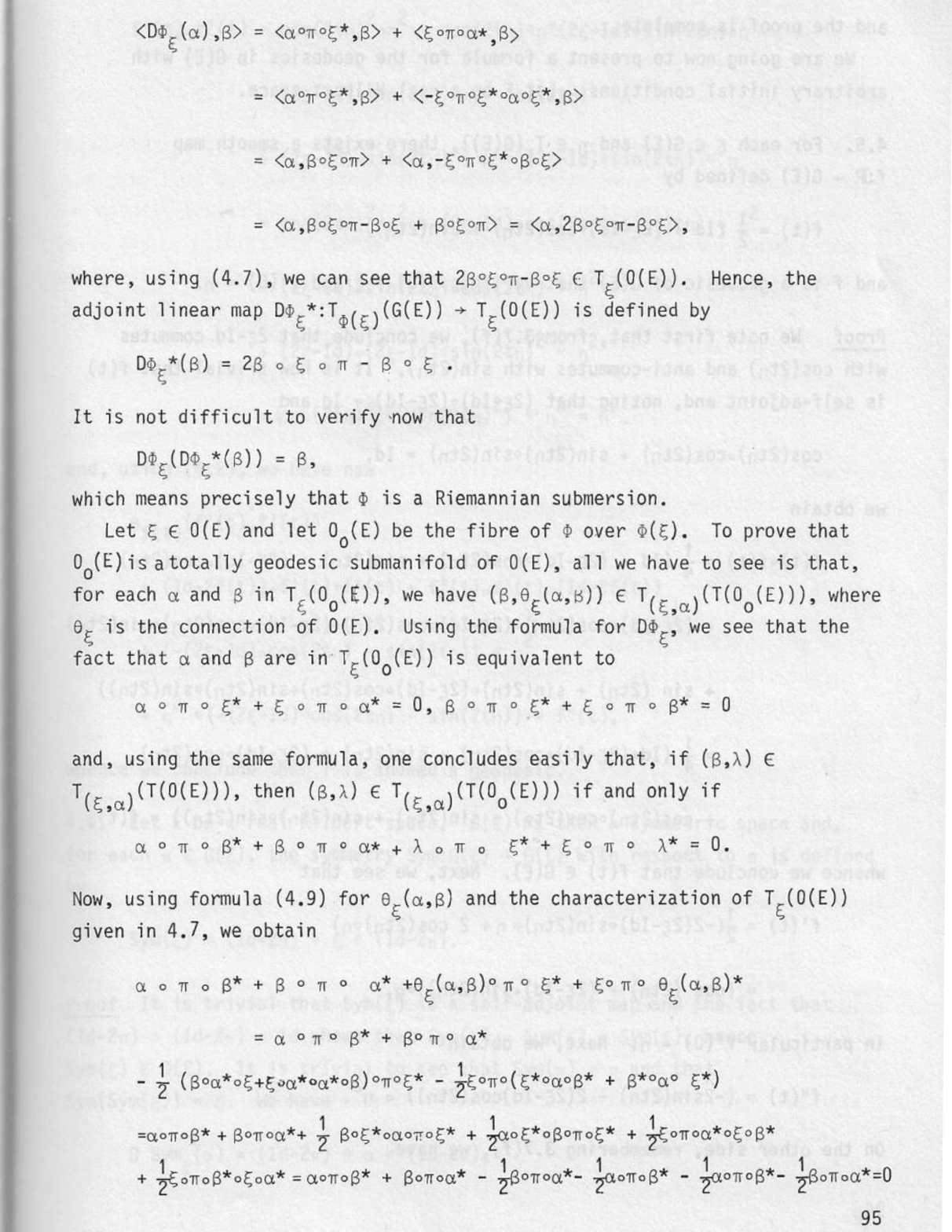}
\vspace*{-7cm}
\clearpage
\thispagestyle{empty}
\mbox{ } \\[-3cm]
\mbox{FIGURE 12. Scanned page 96 of the original paper
published in {\it Res.\ Notes Math.} {\bf131} (1985), 85--102}
\mbox{ } \\[2mm]
\includegraphics[width=18cm]{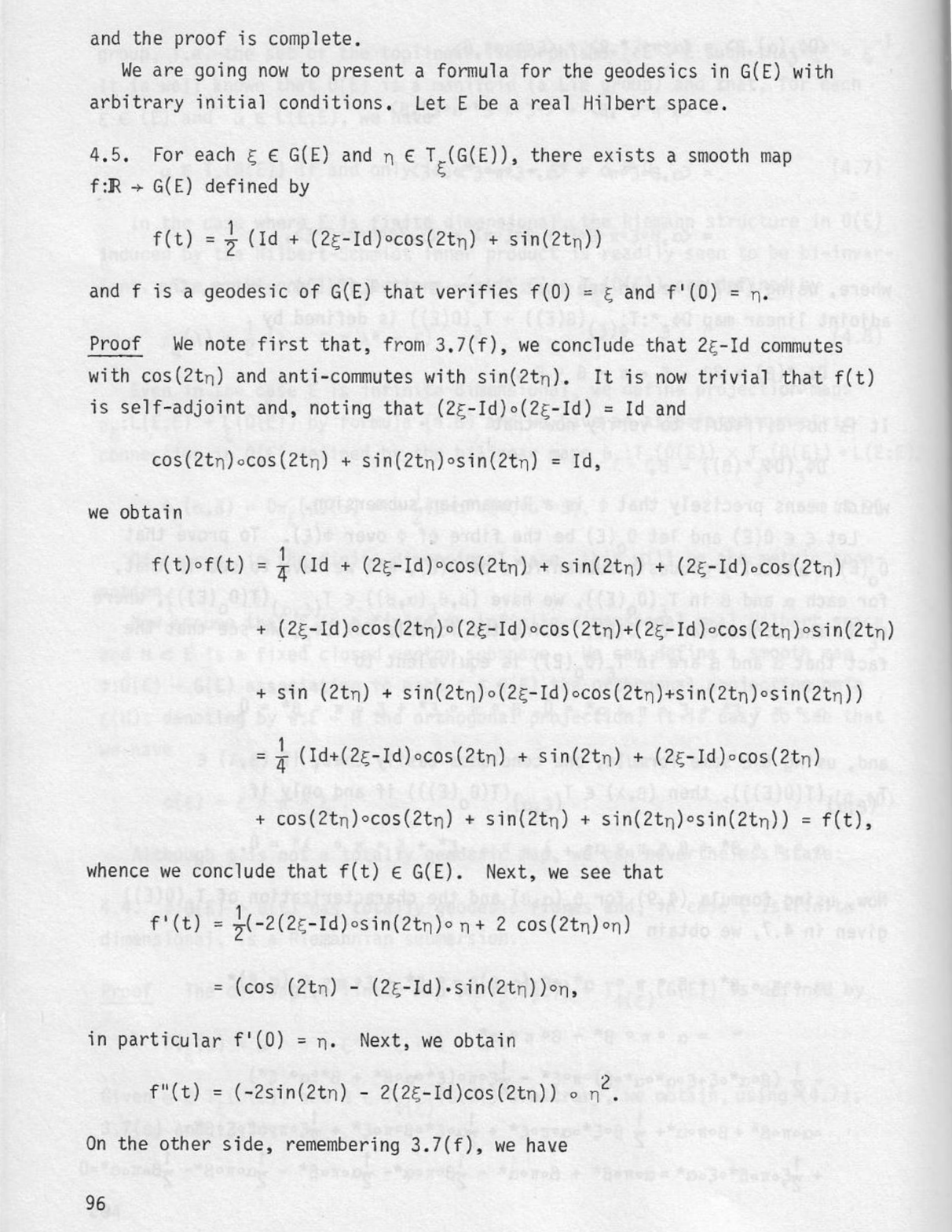}
\vspace*{-7cm}
\clearpage
\thispagestyle{empty}
\mbox{ } \\[-3cm]
\mbox{FIGURE 13. Scanned page 97 of the original paper
published in {\it Res.\ Notes Math.} {\bf131} (1985), 85--102}
\mbox{ } \\[2mm]
\includegraphics[width=18cm]{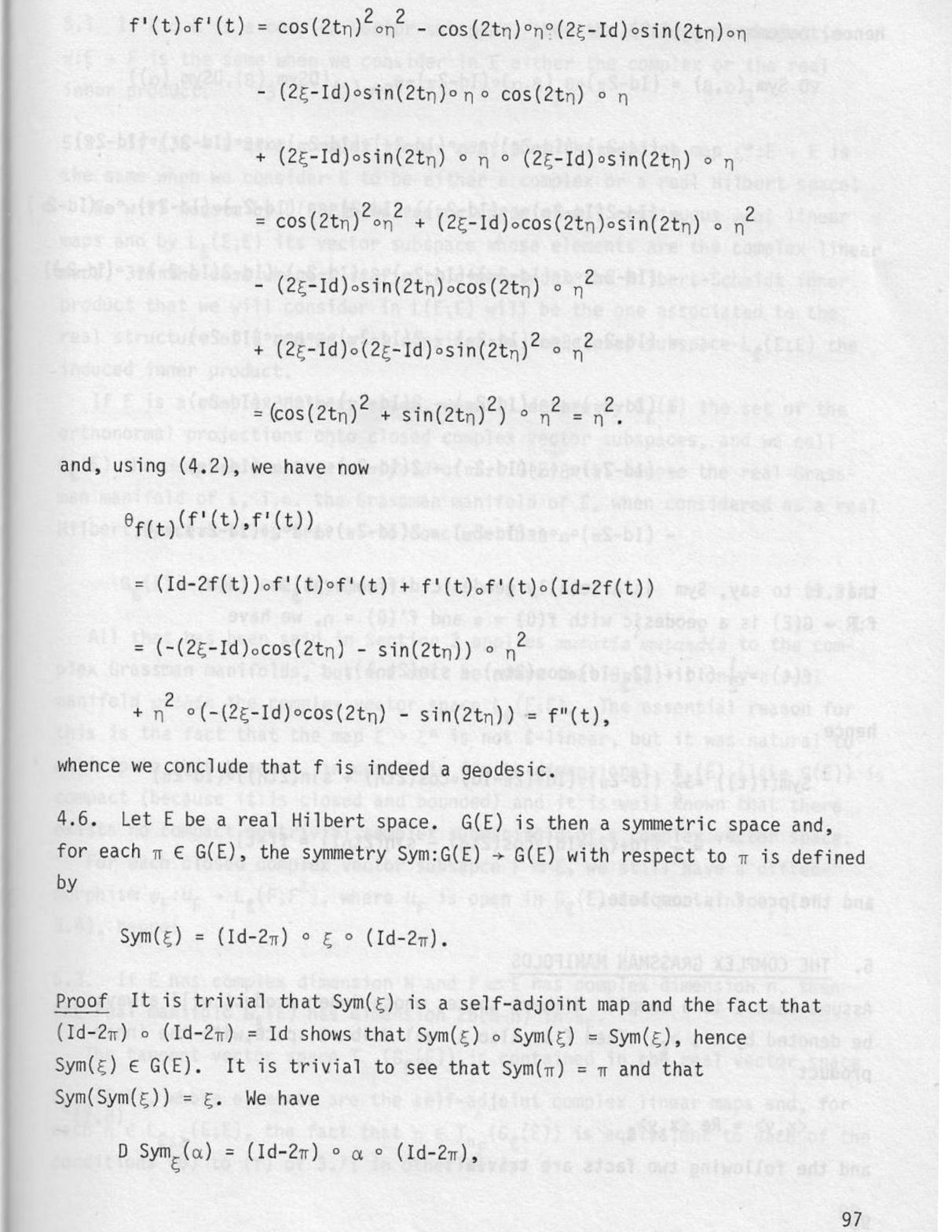}
\vspace*{-7cm}
\clearpage
\thispagestyle{empty}
\mbox{ } \\[-3cm]
\mbox{FIGURE 14. Scanned page 98 of the original paper
published in {\it Res.\ Notes Math.} {\bf131} (1985), 85--102}
\mbox{ } \\[2mm]
\includegraphics[width=18cm]{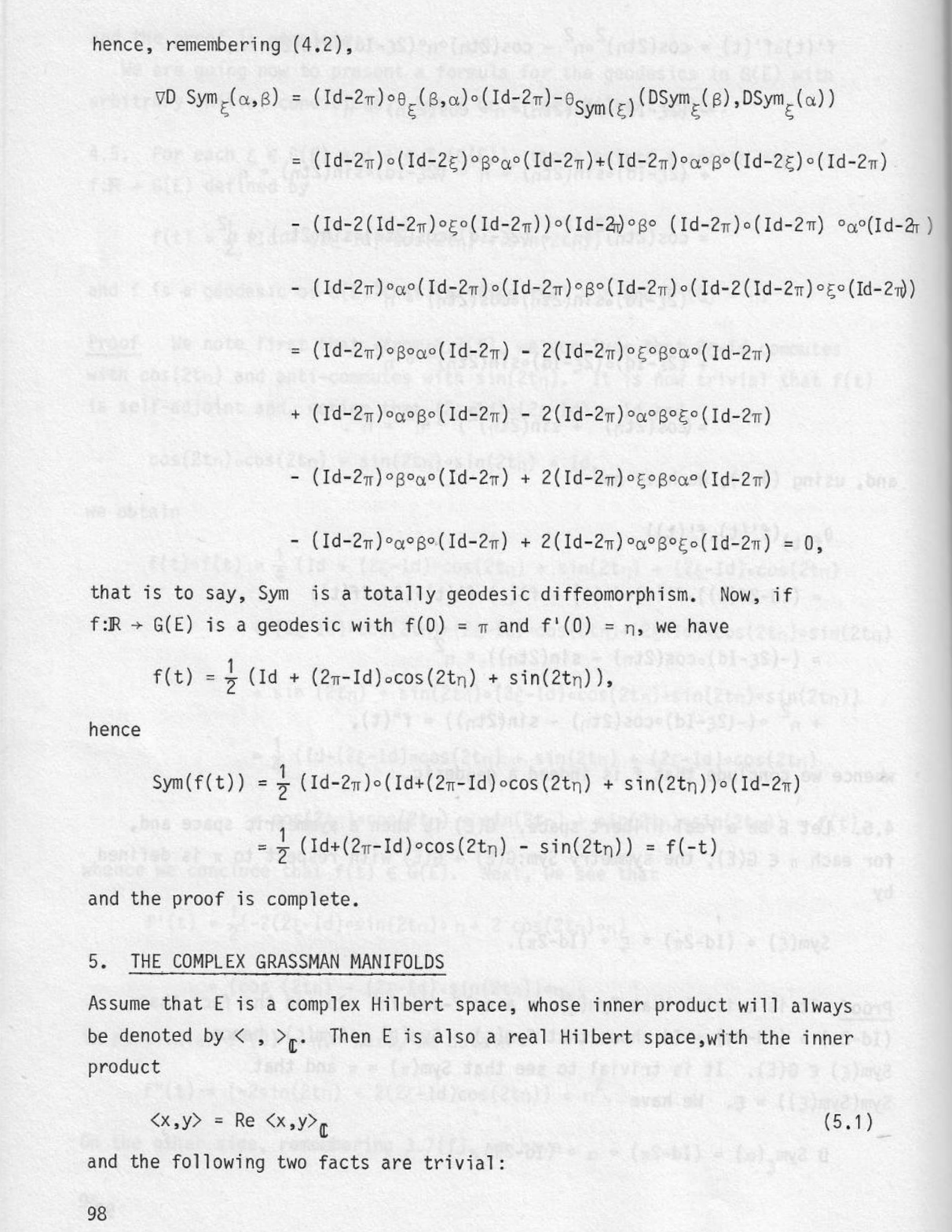}
\vspace*{-7cm}
\clearpage
\thispagestyle{empty}
\mbox{ } \\[-3cm]
\mbox{FIGURE 15. Scanned page 99 of the original paper
published in {\it Res.\ Notes Math.} {\bf131} (1985), 85--102}
\mbox{ } \\[2mm]
\includegraphics[width=18cm]{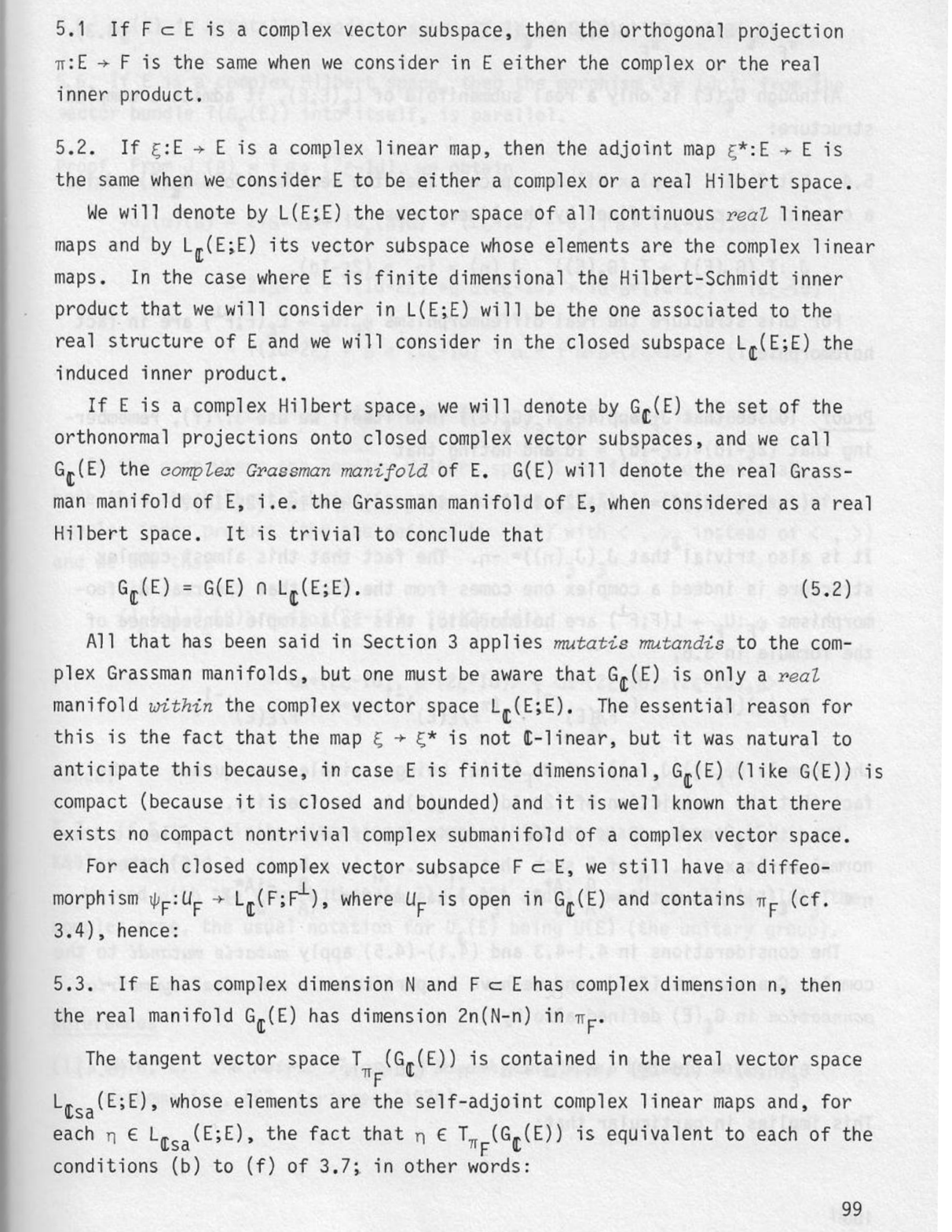}
\vspace*{-7cm}
\clearpage
\thispagestyle{empty}
\mbox{ } \\[-3cm]
\mbox{FIGURE 16. Scanned page 100 of the original paper
published in {\it Res.\ Notes Math.} {\bf131} (1985), 85--102}
\mbox{ } \\[2mm]
\includegraphics[width=18cm]{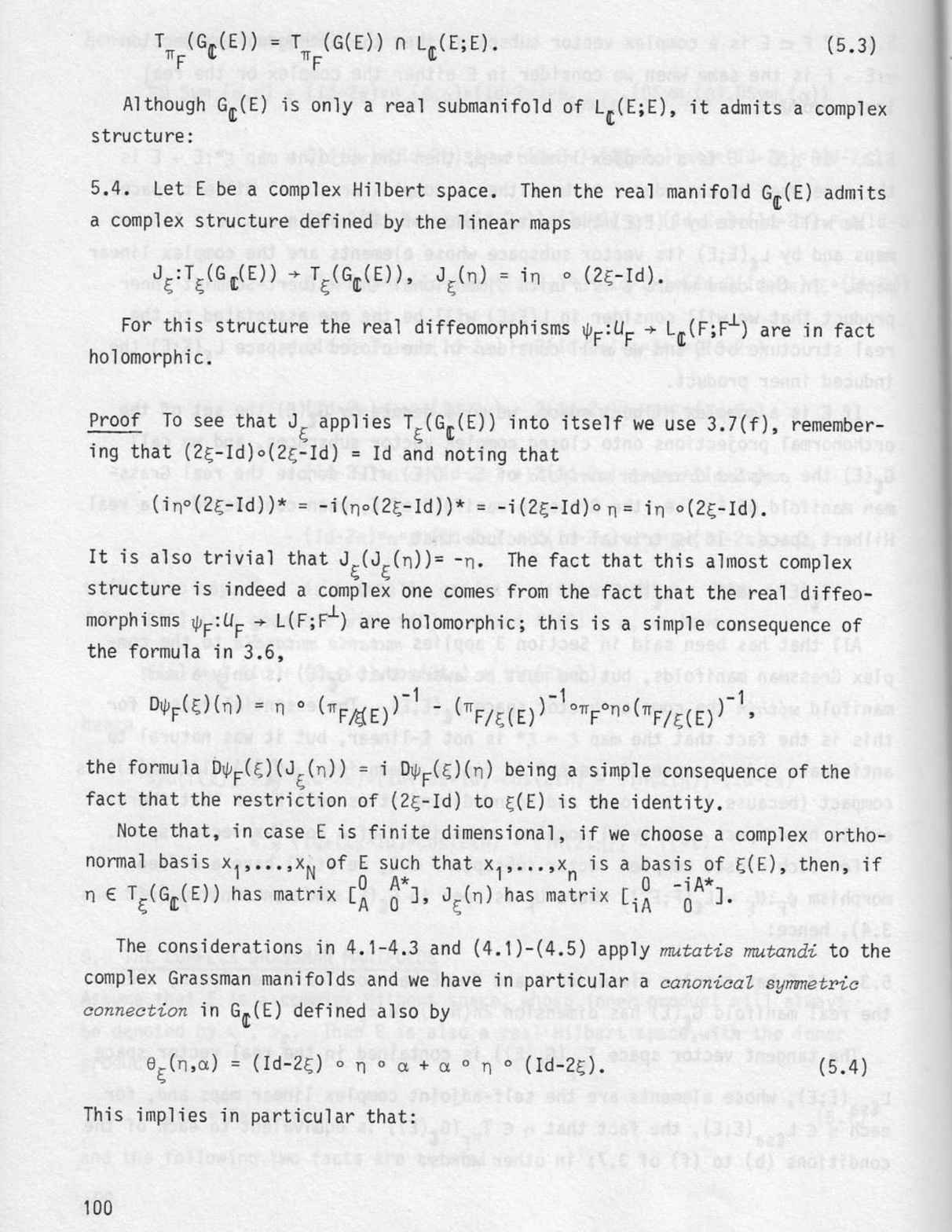}
\vspace*{-7cm}
\clearpage
\thispagestyle{empty}
\mbox{ } \\[-3cm]
\mbox{FIGURE 17. Scanned page 101 of the original paper
published in {\it Res.\ Notes Math.} {\bf131} (1985), 85--102}
\mbox{ } \\[2mm]
\includegraphics[width=18cm]{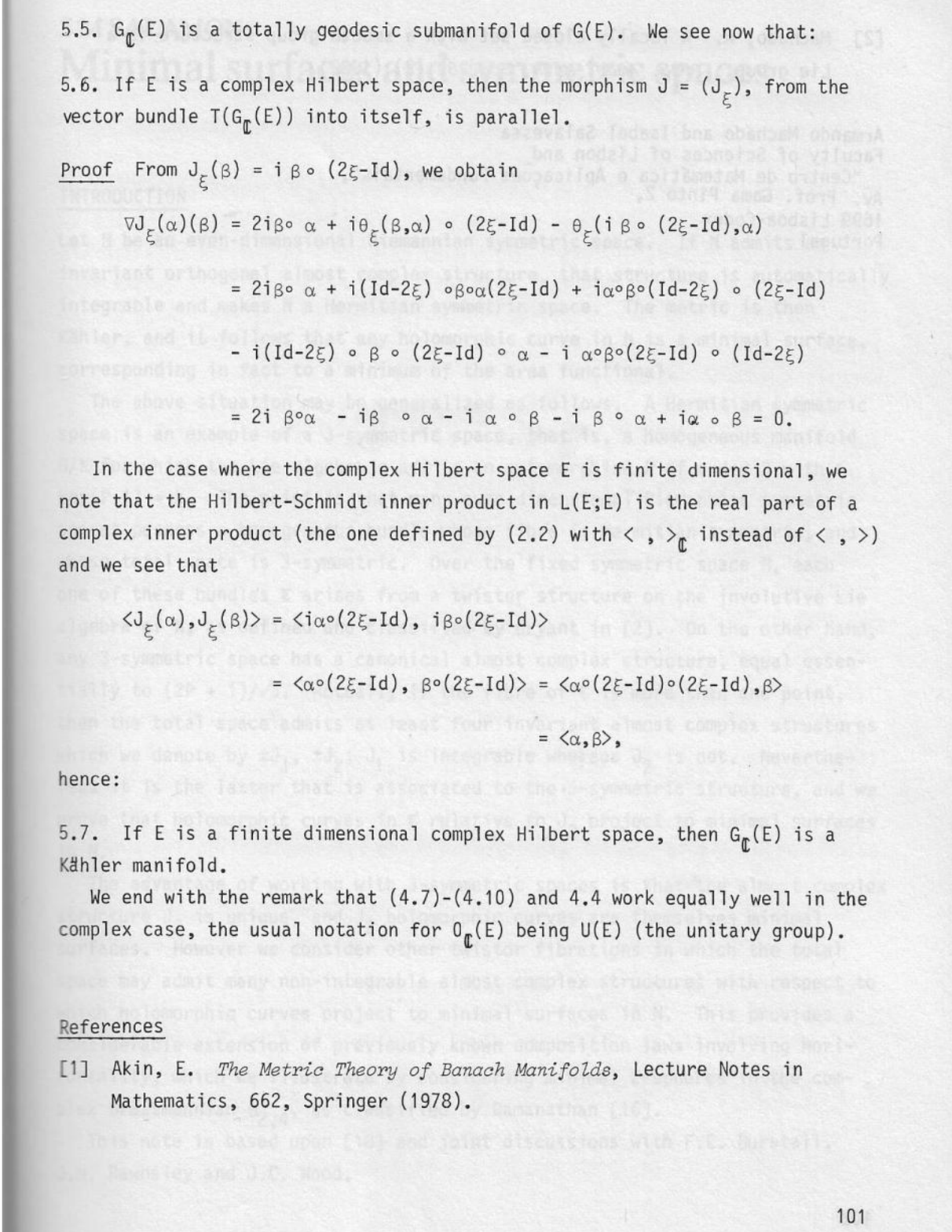}
\vspace*{-7cm}
\clearpage
\thispagestyle{empty}
\mbox{ } \\[-3cm]
\mbox{FIGURE 18. Scanned page 102 of the original paper
published in {\it Res.\ Notes Math.} {\bf131} (1985), 85--102}
\mbox{ } \\[2mm]
\includegraphics[width=18cm]{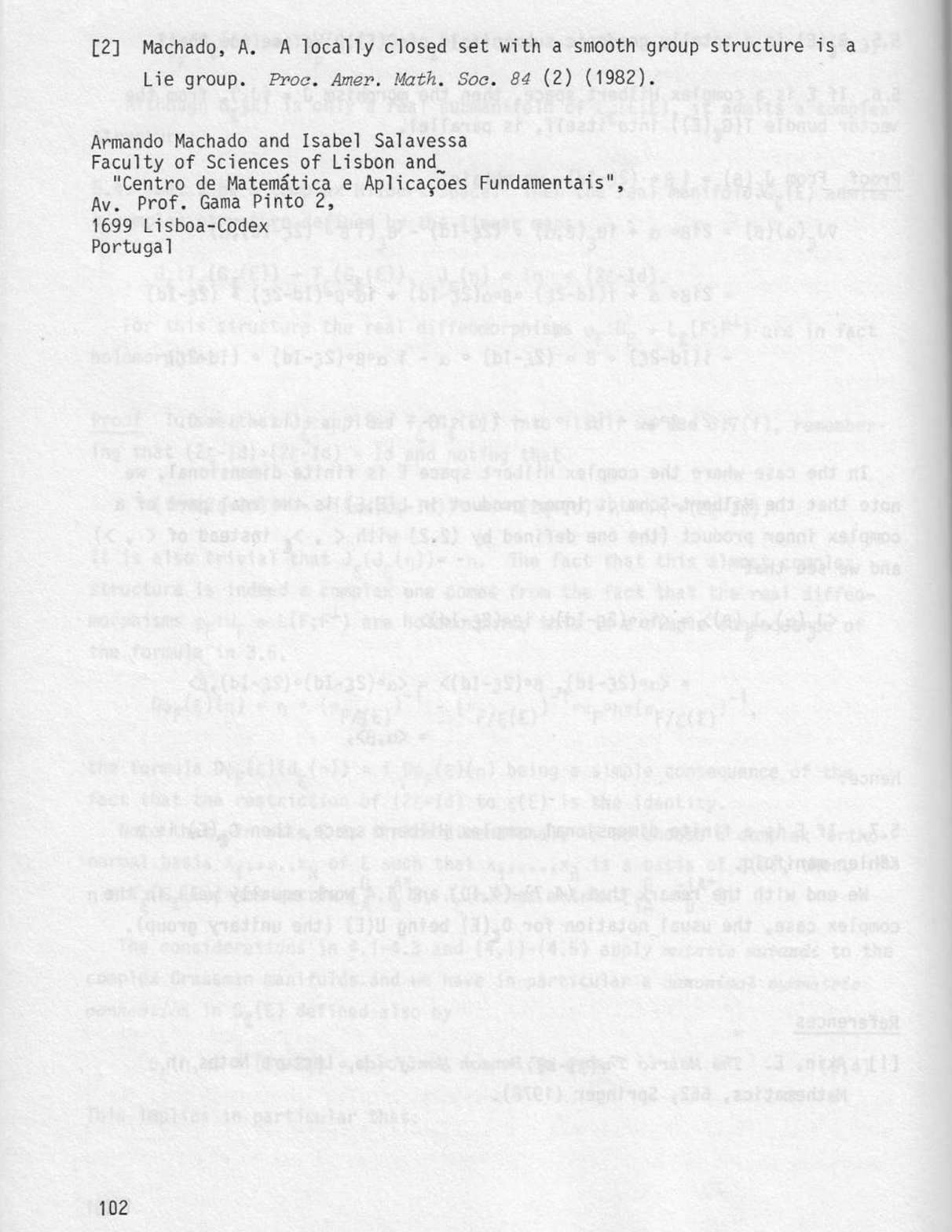}
\vspace*{-7cm}
\end{document}